\documentclass[12pt]{article}
\usepackage[margin=0.75in,top=0.6in]{geometry}
\usepackage[utf8]{inputenc}
\usepackage{amsmath,amsfonts,amssymb}
\usepackage{algorithm}
\usepackage{algpseudocode}
\usepackage{graphicx}
\usepackage{subfig}
\usepackage{float}
\usepackage{multirow}
\usepackage{array}  
\usepackage{booktabs}  
\usepackage{siunitx}  
\usepackage{caption}
\usepackage{hyperref}
\usepackage{nomencl}
\usepackage{multicol}
\makenomenclature
\renewcommand{\nompreamble}{\begin{multicols}{2}}
\renewcommand{\nompostamble}{\end{multicols}}
\date{}
\title{Thermo-Electro-Mechanical Modeling of Power Transmission Line Failures across Four US States}
\begin{document}
\author{
    Prakash KC\thanks{Department of Mechanical Engineering, Michigan State University} \and
    Maryam Naghibolhosseini\thanks{Department of Communicative Sciences and Disorders, Michigan State University} \and
    Mohsen Zayernouri\thanks{Department of Mechanical Engineering \& Department of Statistics and Probability, Michigan State University, Corresponding Author; zayern@msu.edu}
}
\captionsetup{font=small}
\maketitle
\begin{abstract}
    The failure of overhead transmission lines in the United States can lead to significant economic losses and widespread blackouts, affecting the lives of millions. This study focuses on the reliability of transmission lines, specifically examining the effects of wind, ambient temperature, and current demands on lines, incorporating minimal and significant pre-existing damage. We develop a Thermo-Electro-Mechanical Model to analyze the transmission line failures across sensitive and affected states of the United States, integrating historical data on wind and ambient temperature. By combining numerical simulation with historical data analysis, our research assesses the impact of varying environmental conditions on the reliability of transmission lines. Our methodology begins with a deterministic approach to model temperature and damage evolution, using phase-field modeling for fatigue and damage, coupled with electrical and thermal models. Later, we adopt the Probability Collocation Method to investigate the stochastic behavior of the system, enhancing our understanding of uncertainties in model parameters, conducting sensitivity analysis, and estimating the probability of failures over time. This approach allows for a comprehensive analysis of factors affecting transmission line reliability, contributing valuable insights into improving power line's resilience against environmental conditions.  

\end{abstract} 
\noindent \hspace{10mm}\textbf{Keywords:}
{\footnotesize Transmission line, Finite Element Method, Probability Collocation Method, Uncertainty Quantification, Sensitivity Analysis, Probability of Failure}

\nomenclature{\(\theta_{lim}\)}{Threshold Temperature}
\nomenclature{\(\theta_{max}\)}{Maximum Temperature}
\nomenclature{\(\varphi_{lim}\)}{Threshold Damage}
\nomenclature{\(\varphi_{max}\)}{Maximum Damage}
\nomenclature{\(A_{0}\)}{Mean Coefficient}
\nomenclature{\(A_{n}, B_{n}\)}{Frequency Coefficients}
\nomenclature{\(t\)}{Time}
\nomenclature{\(\mathcal{F}\)}{Fatigue Field}
\nomenclature{\(u\)}{Displacement}
\nomenclature{\(Y\)}{Young Modulus}
\nomenclature{\(g_c\)}{Fracture Energy Release Rate}
\nomenclature{\(\gamma\)}{Phase Field Layer Width}
\nomenclature{\(\mathcal{H}'(\varphi), \mathcal{H}_f'(\varphi) \)}{Damage Potentials}
\nomenclature{\(\theta_{c}\)}{Conductor Temperature}
\nomenclature{\(\varphi\)}{Damage}
\nomenclature{\(\theta_{c}\)}{Conductor Temperature}
\nomenclature{\(\theta_{0}\)}{Reference Temperature}
\nomenclature{\(\rho\)}{Material Density}
\nomenclature{\(a\)}{Aging Rate}
\nomenclature{\(q_j\)}{Joule Heating}
\nomenclature{\(\varphi\)}{Damage}
\nomenclature{\(q_c\)}{Convection Cooling}
\nomenclature{\(h\)}{Convective Heat Transfer Coefficient}
\nomenclature{\(v\)}{Wind Velocity}
\nomenclature{\(\rho_{air}\)}{Air Density}
\nomenclature{\(\nu\)}{Kinematic Viscosity of Air}
\nomenclature{\(D\)}{Conductor Diameter}
\nomenclature{\(Re_{D}\)}{Reynolds Number}
\nomenclature{\(Pr\)}{Prandtl Number}
\nomenclature{\(I_b\)}{Base Current}
\nomenclature{\(\theta_b\)}{Base Temperature}
\nomenclature{\(w_b\)}{Base Wind Speed}
\nomenclature{\(q_j\)}{Joule Heating}
\nomenclature{\(\sigma_e\)}{Electrical Conductivity}
\nomenclature{\(\sigma_{E,T}\)}{Non Degraded Electrical Conductivity}
\nomenclature{\(\sigma_0\)}{Electrical Conductivity at Reference Temperature}
\nomenclature{\(P_w\)}{Wind Pressure}
\nomenclature{\(W_w\)}{Wind Load}
\nomenclature{\(C_D\)}{Drag Coefficient}
\nomenclature{\(\alpha\)}{Span Factor}
\nomenclature{\(\theta_w\)}{Wind Direction}
\nomenclature{\(w\)}{Test Function}
\nomenclature{\(A(x)\)}{Cross section Area}
\nomenclature{\(A_0\)}{Undamaged Cross section Area}
\nomenclature{\(J\)}{Current Density}
\nomenclature{\(A_\sigma\)}{Level of Damage}
\nomenclature{\(\delta_t\)}{Time step}
\nomenclature{\(h_B\)}{Bernoulli Random Variable}
\nomenclature{\(p_h\)}{Bernoulli parameter}
\nomenclature{\(H_0\)}{Initial Horizontal Tension}

\printnomenclature

\section{Introduction}

The power grid's components are so interconnected that the failure of a single component can cause a widespread outage. Overhead transmission lines exposed to dynamic weather conditions and current loads are particularly vulnerable. Such vulnerabilities not only lead to frequent blackouts across the United States, impacting millions annually but also disrupt daily life and impose significant economic burdens. A report by the Department of Energy estimated that weather-related power outages cost the US economy between \$18 billion and \$33 billion from 2003 to 2012 \cite{DOE}. Some states such as California, Texas, and Michigan are particularly more vulnerable, with Florida having the highest total number of 25,348,824 affected by power outages\cite{Eaton2019}. Hence, understanding the reliability of transmission cables in such sensitive and affected regions is crucial.

Several studies have investigated the impact of environmental conditions on overhead transmission lines. The effect of thermal stress on the life span of Aluminum Conductor Steel-Reinforced (ACSR) cables was investigated by integrating a strength reliability analysis based on Monte Carlo simulation and fuzzy logic-based ductility analysis\cite{hathout2018impact}. Additionally, the damage caused by wind loading has been extensively studied, affecting the transmission lines and the supporting towers\cite{alminhana2023transmission}. Studies by \cite{dua2015dynamic,stengel2014finite,hamada2011behaviour} have extensively examined the dynamic responses and structural effects of wind on transmission line infrastructure. In \cite{guo2018determination}, the overheating effect of wildfire due to the radiative heat transfer was studied. Similarly, the combined effect of ice and wind was studied in \cite{rossi2020combined}. Other models studied the influence of several weather conditions on the thermal rating of the conductor \cite{bendik2018influence,pytlak2009precipitation,castro2017study}. In addition to environmental factors, the choice of conductor material plays a crucial role in determining the reliability of transmission wires. The most commonly used types include all-aluminum conductor (AAC), all-aluminum alloy conductor (AAAC), and aluminum conductor steel-reinforced (ACSR) \cite{zainuddin2020review}. 

These studies laid a solid basis for understanding the failure mechanisms of power transmission. While much of this research has focused on the thermal and electrical modeling of transmission lines, none considered the mechanical model that relates to fatigue and damage that influence the overall conductivity of the material. Transmission lines are subjected to dynamic wind, temperature, and electrical loads, which can induce fatigue and damage to the wires. This damage compromises the conductivity of the conductor, leading to increased heat generation. Furthermore, the cyclic nature of the loading can accelerate damage, ultimately precipitating the early failure of the wire. The development of models that integrate the electrical and thermal properties of overhead conductors with advanced formulations in damage and fatigue mechanics can significantly enhance the predictability of transmission line failures during long-term operations.

Phase fields have emerged as a powerful research area in damage and fatigue modeling. Originally formulated to address fluid separation problems \cite{cahn1958free}, phase-field models exploit their capability to represent sharp interfaces through a smooth, continuous field, broadening their application to various multiphase problems. Over the years, these models have been crucial in simulating brittle \cite{tanne2018crack,li2015phase,borden2014higher,miehe2010phase,miehe2010thermodynamically,de2021data}, ductile \cite{seiler2020efficient,miehe2016phase,kuhn2016phase,ambati2016phase,ambati2015phase}, isothermal fatigue fracture \cite{barros2021integrated}, and non-isothermal fatigue fracture mechanisms\cite{boldrini2016non}, and accurately capturing phenomena like crack initiation, propagation, branching, and coalescence—events commonly observed in dynamic fractures. Recently, phase field models have been integrated with electrical and thermal models to study various phenomena, including the breakdown of polymers under alternating voltages \cite{zhu2020phase}, transitions from paraelectric to ferroelectric states \cite{woldman2019thermo}, and the breakdown of polymer-based dielectrics \cite{shen2019phase}. Based on these applications, incorporating the phase field model with electrical and thermal models becomes crucial in accurately predicting the reliability of transmission lines.  Further incorporation of stochastic analysis to account for variables such as weather conditions and initial damage in practice can enhance the precision of lifespan predictions. This approach enables a more comprehensive understanding of the factors that affect transmission line durability in real-world scenarios.

The study of transmission line reliability requires an understanding of the interplay between mechanical, thermal, and electrical components, along with the uncertainties. While the phase field model is accepted for damage evolution modeling, dislocation dynamics can capture the detail of the microstructural mechanism \cite{de2021atomistic}. Studies by \cite{chhetri2023comparative, de2023machine} shown that material behavior can be significantly affected by dislocation interaction resulting in failure under dynamic conditions. The material behavior under the influence of high temperature deviates from normal behavior to more complex behavior that includes viscous and elastic behavior. Fractional visco-elasto-plastic models \cite{suzuki2022general} can capture these behaviors resulting in better structural analysis \cite{suzuki2016fractional}, accurate damage evolution \cite{suzuki2021thermodynamically}, and large-scale behavior \cite{suzuki2023fractional}. Integrating these fractional-order models can efficiently improve the failure prediction of the transmission line. Further, incorporating fractional order as a random variable can better capture the system response\cite{kharazmi2019operator}.

The classical Monte Carlo method \cite{fishman2013monte,smith2013uncertainty} is a standard benchmark for studying stochastic solutions that compute the Quality of Interest (QoI) in a straightforward fashion. However, this method is characterized by a slow convergence rate, requiring a large number of realizations and thus making the process computationally expensive. This limitation highlights the need for more efficient computational approaches in stochastic analysis. Other established methods like Polynomial Chaos \cite{xiu2002modeling,xiu2002wiener}, or its generalization via Galerkin projection \cite{stefanou2009stochastic,babuvska2005solving,babuska2004galerkin}, require modifications to the governing equations for stochastic analysis. This requirement makes these methods intrusive and potentially impractical for complex problems since such modifications may not be feasible or could overly complicate the analysis. To address this, non-intrusive techniques must be used; one such technique is the Probabilistic Collocation Method (PCM) \cite{babuvska2007stochastic,xiu2005high}. PCM preserves the simplicity of solution structures and allows for independently sampled realizations, thereby achieving better convergence compared to the traditional Monte Carlo method. Although the curse of dimensionality challenges PCM due to tensorial products, this issue can be effectively addressed using techniques such as sparse grids \cite{smolyak1963quadrature} or active subspace methods \cite{constantine2017global,constantine2015exploiting,constantine2014active}.

In a recent study by our group \cite{demoraes2024thermoelectromechanicalmodellongtermreliability}, the authors initially examined a deterministic solution and later implemented the Probabilistic Collocation Method (PCM) for stochastic analysis. Our integrated approach incorporated a phase field model with electrical and thermal analyses to predict the lifespan of overhead transmission lines. However, a more comprehensive analysis involving historical data across various specific states would provide a more realistic assessment. The reliability and durability of transmission lines vary under different weather conditions, highlighting the need to adjust the weather condition assessments to particular cases. Furthermore, analyzing the lifespan of transmission lines with insignificant initial damage could offer valuable insights into the extent to which pre-existing damage influences overall reliability.
In this study, we incorporate several realistic scenarios for a detailed historical analysis of wind and temperature data by focusing on four specific U.S. states: Texas, California, Michigan, and Florida. Additionally, we assess the system’s reliability under minimal damage to establish a baseline and examine how precursor damage affects the reliability of transmission lines.

This paper is organized as follows: In Section 2, we present the problem statement including four representative scenarios along with the historical data analysis for wind and temperature of each specific scenario. We present a thermo-electro-mechanical model for the failure of transmission lines and discuss each model in section 3. Then, we discuss the deterministic solution through the finite-element method in Section 4. The stochastic methods are discussed in Section 5, where we present the PCM building blocks for the uncertainty, sensitivity, and probability of failure analyses. Finally, we address the concluding remarks in Section 6.

\section{Problem Statement}
\label{sec:problem statement}

The operating temperature and damage parameters are crucial for the reliable operation of overhead transmission lines. High temperatures can cause annealing and sagging of the wire, while cyclic loading in the presence of initial damage on the wire leads to crack initiation and propagation, ultimately affecting the wire's life and reliability. These issues often remain undetected until a cable rupture occurs. Additionally, environmental conditions including ambient air temperature, wind speed, and current demand can adversely affect these parameters, pushing them beyond acceptable levels. 

Therefore, operating temperature and damage are considered the primary factors for our analysis. Along with physical effects, material parameters, and loading conditions, these elements determine the failure state of the transmission line over its lifespan. This study aims to understand the multiphysics effects on these critical factors, which are crucial for assessing the transmission line's reliability. 

\subsection{Representative Scenarios}

Previous research on the reliability of transmission lines often relied on parameterized values for wind and temperature, which lack precision. This study aims to improve this by incorporating real data from four representative scenarios: Texas, California, Michigan, and Florida. We source wind and temperature data from \cite{noaa_ccd_2024} and \cite{nws_2024}. The data for each state are presented in Table \ref{tab:wind_temp_data}.
\begin{flushleft}
\begin{table}[H]
\caption{Wind speed (ft/s)\cite{noaa_ccd_2024} and Temperature (K) \cite{nws_2024} Data by State.}
\renewcommand{\arraystretch}{2.5} 
\small
\begin{tabular}{cc*{12}{p{0.8cm}}} 
\hline
\multirow{2}{*}{Texas} & ft/s & 11.59 & 12.32 & 12.76 & 13.05 & 12.17 & 10.41 & 9.09 & 8.65 & 9.24 & 10.12 & 10.85 & 11.15 \\
\cline{2-14}
                       & K    & 287.54 & 286.98 & 294.82 & 294.71 & 298.71 & 301.21 & 303.48 & 303.48 & 300.37 & 295.21 & 292.04 & 286.15 \\
\hline
\multirow{2}{*}{California}  & ft/s & 2.05 & 2.79 & 3.08 & 3.52 & 3.23 & 3.08 & 3.08 & 2.93 & 2.49 & 2.20 & 2.05 & 2.05 \\
\cline{2-14}
                       & K    & 286.65 & 288.04 & 287.76 & 290.09 & 293.37 & 294.21 & 295.71 & 297.54 & 297.15 & 295.15 & 289.76 & 287.65 \\
\hline
\multirow{2}{*}{Michigan} & ft/s & 16.13 & 15.40 & 15.25 & 15.25 & 13.35 & 12.32 & 11.59 & 10.71 & 11.29 & 13.49 & 15.11 & 15.40 \\
\cline{2-14}
                          & K    & 273.37 & 272.09 & 278.26 & 280.98 & 287.43 & 294.87 & 298.43 & 296.21 & 290.59 & 283.71 & 280.76 & 273.65 \\
\hline
\multirow{2}{*}{Florida} & ft/s & 12.61 & 12.91 & 14.23 & 14.23 & 13.20 & 11.00 & 10.56 & 10.41 & 10.85 & 12.61 & 12.91 & 12.17 \\
\cline{2-14}
                         & K    & 295.09 & 295.98 & 298.26 & 300.87 & 299.48 & 302.26 & 303.09 & 302.76 & 301.93 & 301.04 & 298.43 & 293.93 \\
\hline
\end{tabular}
\label{tab:wind_temp_data}
\end{table}
\end{flushleft}

\subsection{Discrete Fourier Transform}

\begin{figure}[H]
	\centering
    {\includegraphics[width=0.24\textwidth]{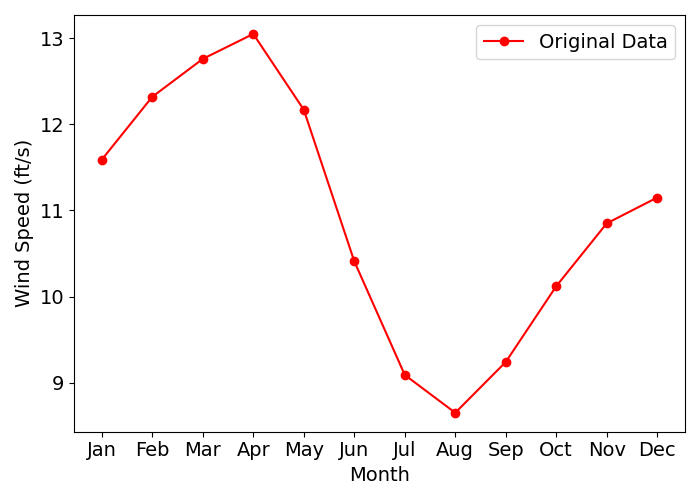}}
    {\includegraphics[width=0.24\textwidth]{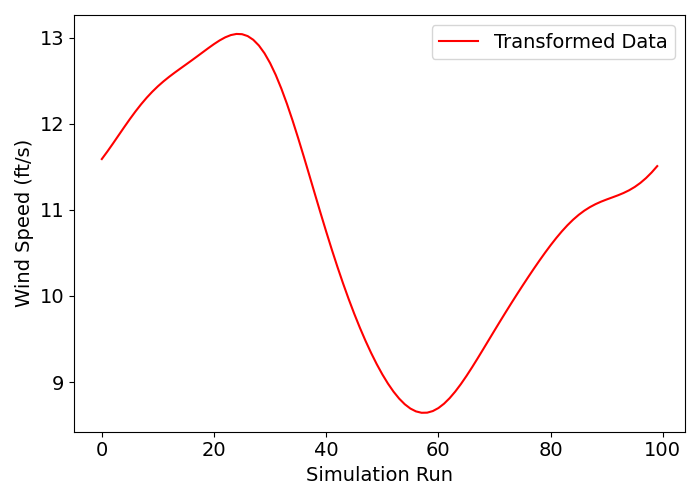}}
    {\includegraphics[width=0.24\textwidth]{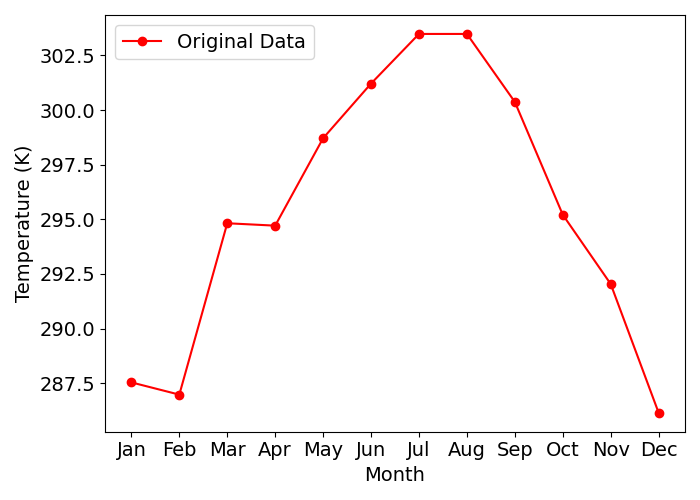}}
    {\includegraphics[width=0.24\textwidth]{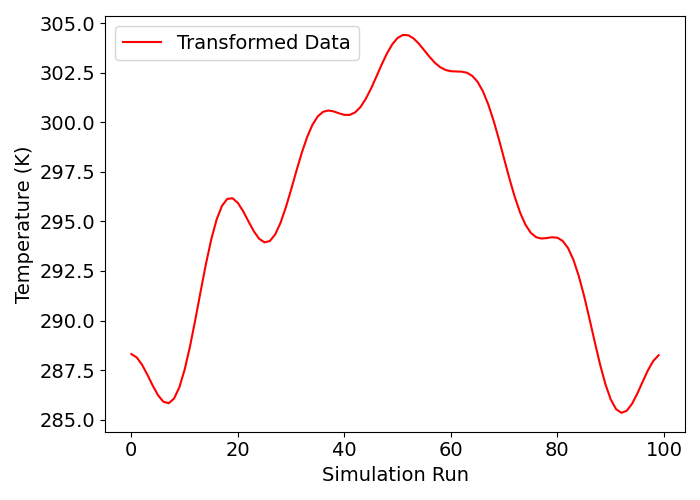}}
	\caption{Discrete and continuous variation of wind and temperature data of Texas}
	\label{fig:Texas}
\end{figure}

\begin{figure}[H]
	\centering
    {\includegraphics[width=0.24\textwidth]{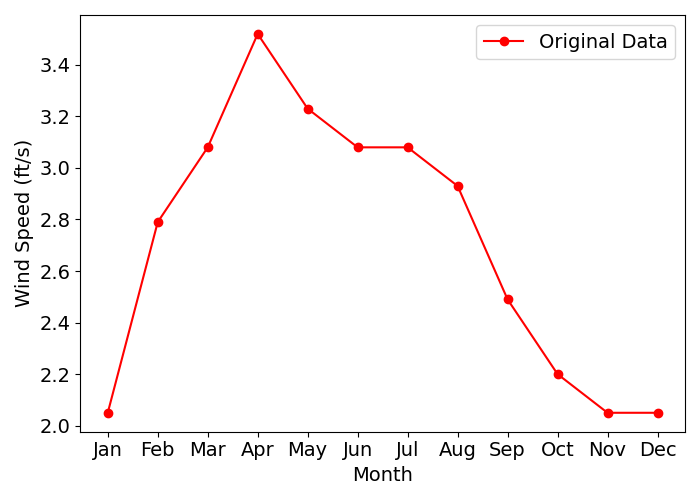}}
    {\includegraphics[width=0.24\textwidth]{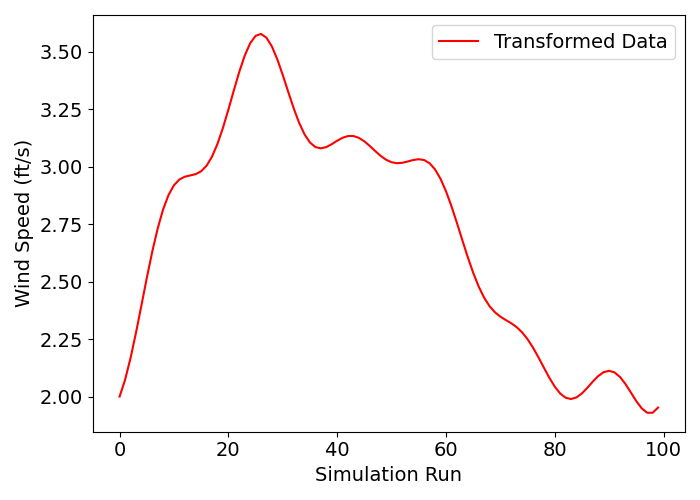}}
    {\includegraphics[width=0.24\textwidth]{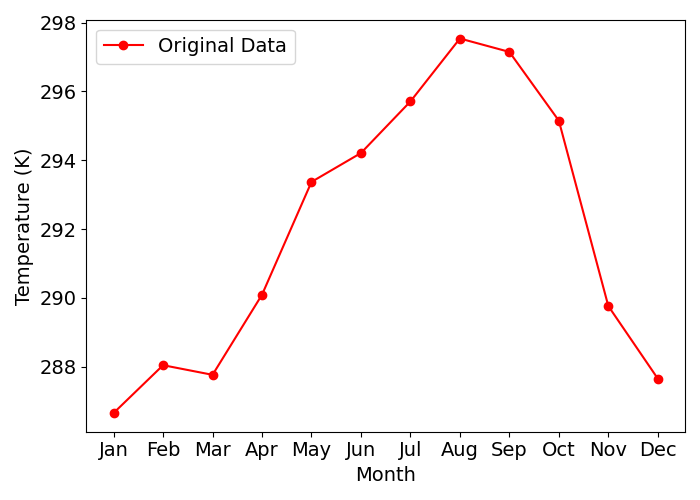}}
    {\includegraphics[width=0.24\textwidth]{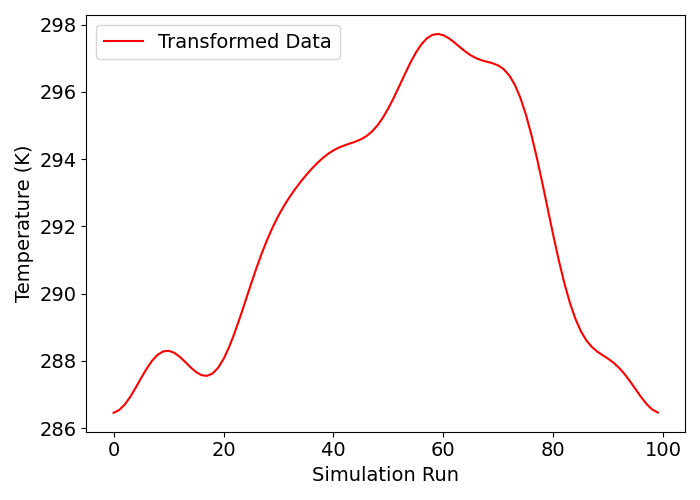}}
	\caption{Discrete and continuous variation of wind and temperature data of California}
	\label{fig:California}
\end{figure}

\begin{figure}[H]
	\centering
    {\includegraphics[width=0.24\textwidth]{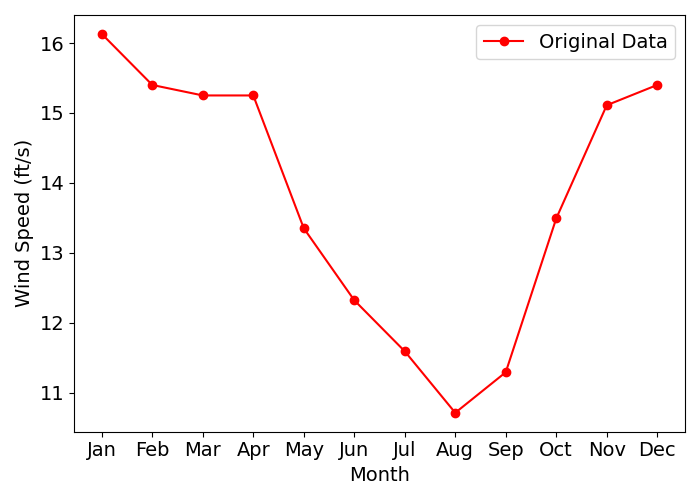}}
    {\includegraphics[width=0.24\textwidth]{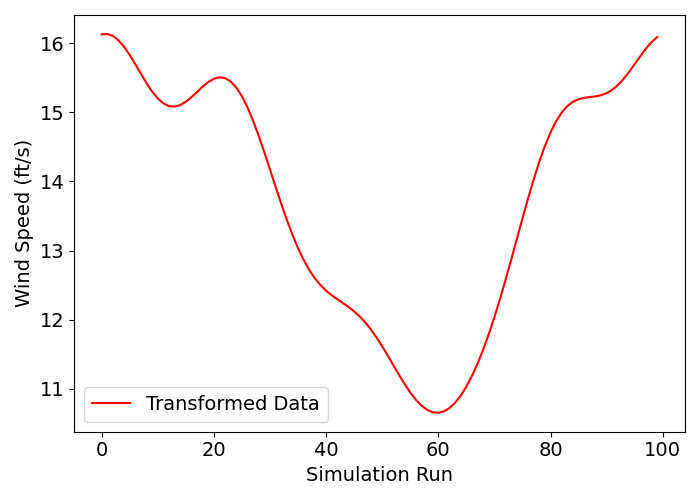}}
    {\includegraphics[width=0.24\textwidth]{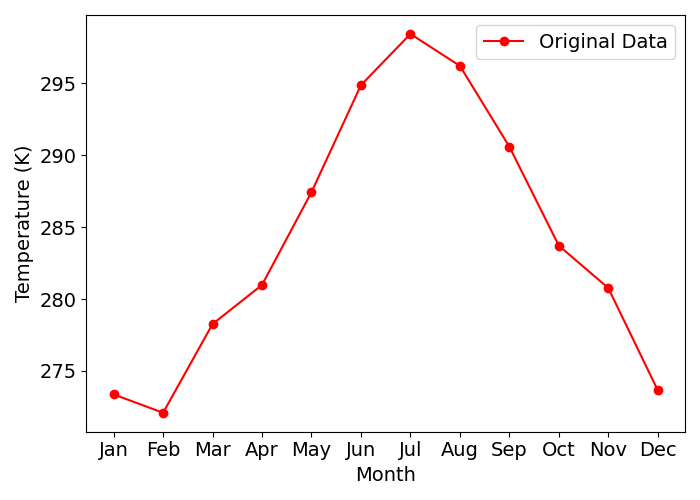}}
    {\includegraphics[width=0.24\textwidth]{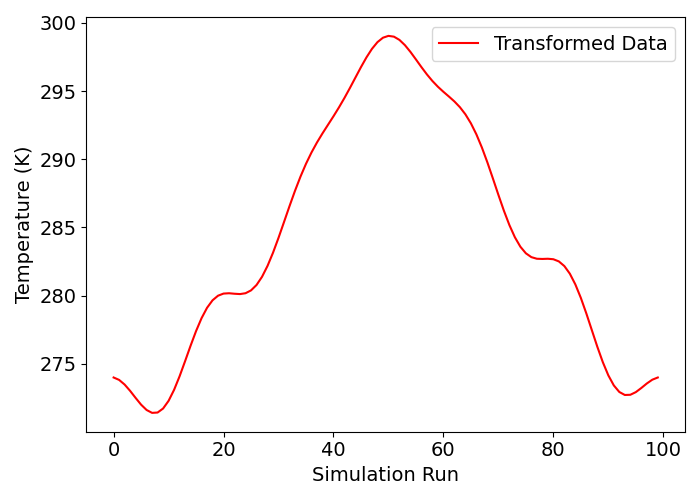}}
	\caption{Discrete and continuous variation of wind and temperature data of Michigan}
	\label{fig:Michigan}
\end{figure}

\begin{figure}[H]
	\centering
    {\includegraphics[width=0.24\textwidth]{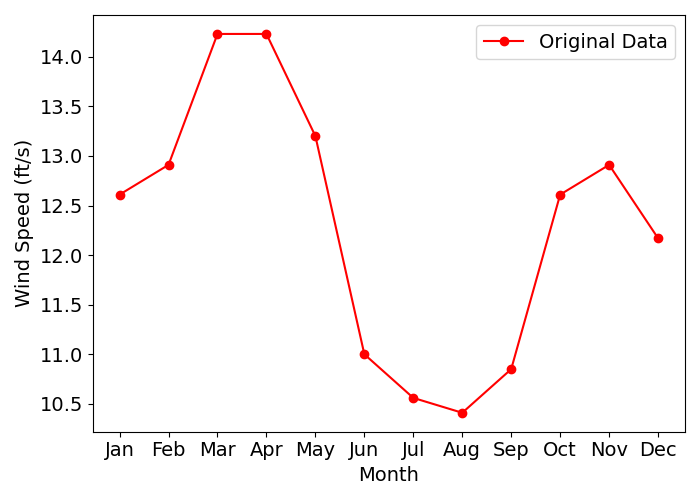}}
    {\includegraphics[width=0.24\textwidth]{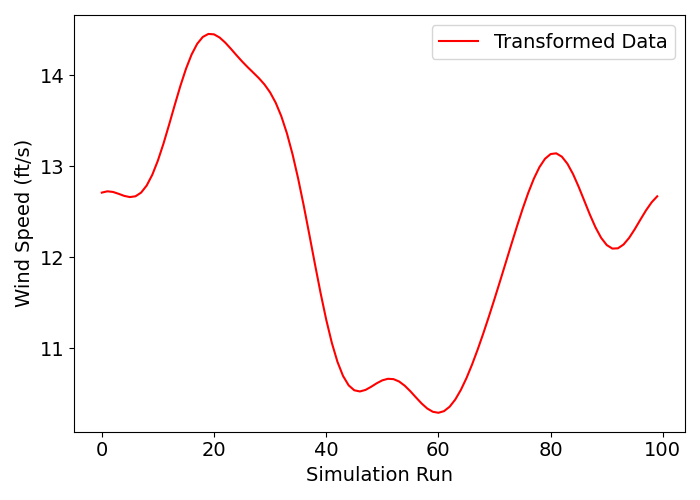}}
    {\includegraphics[width=0.24\textwidth]{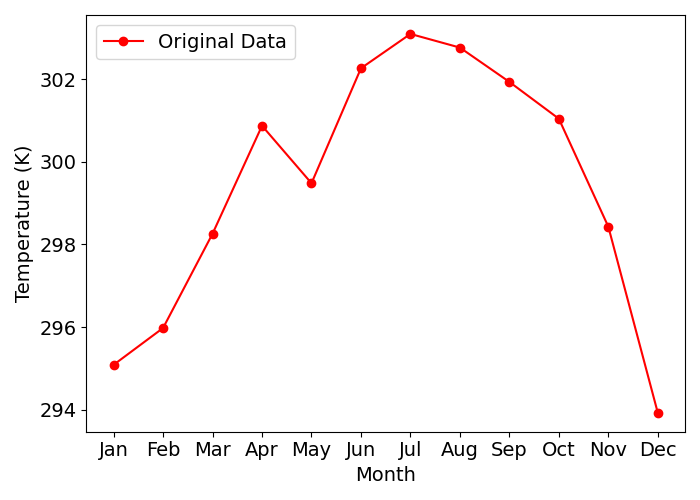}}
    {\includegraphics[width=0.24\textwidth]{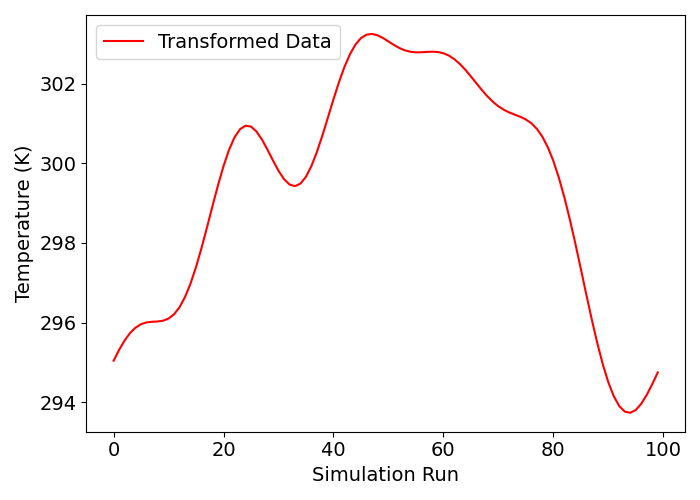}}
	\caption{Discrete and continuous variation of wind and temperature data of Florida}
	\label{fig:Florida}
\end{figure}

The available data for wind and temperature are discrete. To obtain the loading conditions, we implement the Discrete Fourier Transform (DFT) and Fourier series. The DFT is defined by the following expression, 
\begin{equation}
    X_k = \sum_{n=0}^{N-1} x_n e^{-i2\pi \frac{k}{N}n}
\end{equation}
where $x_n$ is the $n$th sample of $x$ and $N$ is the total number of sample. We obtain the mean, $A_0$, and frequency coefficients, $A_n$ and $B_n$ for each sample from DFT. The mean and frequency coefficients are used in the Fourier series to obtain the cyclic loading equation using:
\begin{equation}
    f(t) = A_0 + \sum_{n=1}^{N/2} [A_n \cos(2\pi n \frac{t}{T}) + B_n \sin(2\pi n \frac{t}{T})]
\end{equation}
where $t$ is scaled to match the original time domain, and $T$ represents the period of the data set. The time is structured such that 100 iterations correspond to the cyclic loading condition of 12 months which is later implemented in our thermo-electro-mechanical model for failure analysis. The discrete (original data) and continuous data (obtained using the Fourier series analysis) plots of each state are shown in Figure (\ref{fig:Texas}),(\ref{fig:California}),(\ref{fig:Michigan}), and (\ref{fig:Florida}) 

\subsection{Reliability of Transmission Lines}

In this study, we focused on assessing the reliability of transmission lines, using a reliability model that revolves around a limit state function $g(R, S; t)$, where $R$ represents a threshold, $S$ denotes the maximum value of a dependent variable, and $t$ is time. Our analysis specifically addresses the reliability concerning the maximum temperature or maximum damage a transmission line can withstand before failure. We define $S$ as the maximum temperature or maximum damage of the transmission line, $\theta_{\text{max}}$ or $\varphi_{\text{max}}$ respectively, at any given time $t$, and $R$ as the critical temperature or damage limit, $\theta_{\text{lim}}$ or $\varphi_{\text{lim}}$, beyond which the integrity of the transmission line is compromised. Thus, the limit state function is defined as:

\begin{equation}
g(\theta_{\text{lim}}, \theta_{\text{max}}; t) = \theta_{\text{lim}} - \theta_{\text{max}}(t)
\end{equation}
\begin{equation}
g(\varphi_{\text{lim}}, \varphi_{\text{max}}; t) = \varphi_{\text{lim}} - \varphi_{\text{max}}(t)
\end{equation}

Generally, both $R$ and $S$ are treated as random variables in reliability analysis, however, in this study, the threshold limit is held constant. Variability is completely attributed to $\theta_{\text{max}}$ or $\varphi_{\text{max}}$, which are functions of time and space along the domain. Thus, the probability of failure $p_f(t)$, can be defined as:

\begin{equation}
p_f(t) = P\{g(\theta_{\text{lim}}, \theta_{\text{max}}; t) < 0\}.
\end{equation}
\begin{equation}
p_f(t) = P\{g(\varphi_{\text{lim}}, \varphi_{\text{max}}; t) < 0\}.
\end{equation}

The thermo-electro-mechanical model acts as a computational tool to evaluate the operating temperature and damage of the transmission line over time. Using the model in the stochastic framework, we obtain the probability of failure of the transmission line.

\section{Methodology}

\begin{figure}[H]
	\centering
	\includegraphics[width=0.45\textwidth]{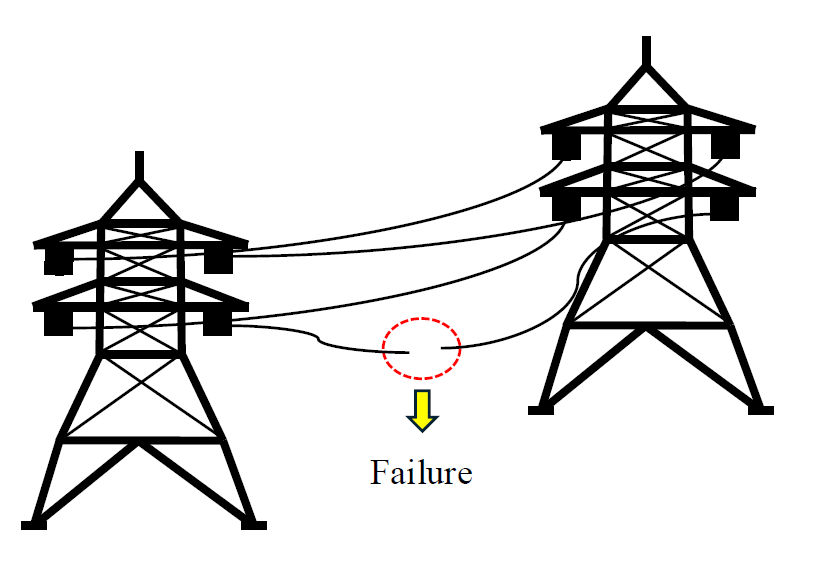}
	\caption{\textit{Schematic representation of Transmission lines.}}
	\label{fig:tower}
\end{figure}

Figure \ref{fig:tower} illustrates the transmission line mounted on the transmission tower, highlighting the cable failure. 

\begin{figure}[H]
	\centering
	\includegraphics[width=0.8\textwidth]{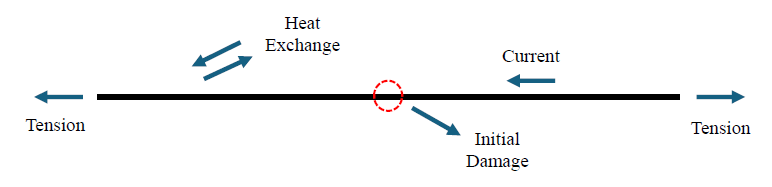}
	\caption{\textit{One-dimensional representation of transmission line.}}
	\label{fig:line}
\end{figure}

In our study, we examine a specific segment of transmission cable. The segment is supported by two transmission towers and the sag is maintained within permissible limits. Mechanically, we treat the cable as a one-dimensional body under tension, using the projected span as its effective operational length for simplification as shown in Figure \ref{fig:line}. While sag can be an essential consideration, we only account for its influence on horizontal tension due to temperature variations. We do not explicitly model the sag; instead, we focus on understanding how the horizontal tension contributes to material damage and fatigue over time.

The horizontal tension in the cable is primarily influenced by its operating temperature. As the temperature increases, the cable experiences thermal elongation, which increases the sag and consequently reduces the horizontal tension. Conversely, when the cable cools, it contracts, resulting in increased tension. Additionally, this repeated elongation and contraction due to sagging can lead to fatigue in the conductor, accelerating material damage and eventually leading to mechanical failure.

The temperature of the conductor is driven by Joule heating and convective cooling from the wind acting on it. Additionally, any damage to the conductor can increase its resistivity, further increasing its temperature. The increase in temperature accelerates the degradation and aging of the material. This illustrates how interconnected temperature and damage are affecting the longevity of power lines.

The electrical model, here, does not consider the electrical design of transmission lines. Rather, it focuses on the effect of the electric current passing through a single cable on the material properties, independent of the voltage levels delivered to consumers. The primary consideration of this approach is that resistive losses along the line increase the material's temperature and result in a voltage drop from one end of the line to the other.

In this study, we propose that any damage to the material increases its resistivity, further increasing the temperature-induced resistivity effects and leading to greater heat generation due to Joule heating. The heat generation elevates the temperature of the material and accelerates aging. This creates a positive feedback loop where damage accelerates heating, which in turn further increases the temperature and accelerates damage, potentially leading to early failure due to thermal runaway or mechanical degradation. To model these cycles over time, we use a quasi-static approach, assuming the system reaches equilibrium more quickly than the change in loading conditions.

\subsection{Mechanical Model}

We implement a non-isothermal phase-field framework to model damage and fatigue based on principles in \cite{boldrini2016non}. This framework consists of two partial differential equations (PDEs) which govern the evolution of displacement, denoted by $u$, and damage, denoted by $\varphi$, along with an ordinary differential equation (ODE) for the fatigue field, $\mathcal{F}$.
The damage phase field represents the volumetric fraction of degraded material, with $\varphi = 0$ indicating undamaged material and $\varphi = 1$ indicating complete fracture, while values $0 < \varphi < 1$ represent the intermediate level of damage. Damage follows an Allen-Cahn-type equation, obtained along with the equilibrium equation for $u$ by applying the principle of virtual power and ensuring thermodynamic consistency through entropy inequalities. 

The fatigue field, $\mathcal{F}$, is considered an internal variable. Its evolution is described by an equation derived from constitutive relations that must satisfy the entropy inequality to ensure thermodynamic consistency for all admissible processes.

We adopt a 1-D model for the mechanical body, which occupies the domain $\Omega \in \mathbb{R}$ at time $t \in (0, T]$. From the governing equations, specific material evolution behaviors can be derived by selecting particular free-energy potentials. An alternative way is considering the free-energy function:

\begin{equation}
\label{eq:free_energy}
\Psi (\nabla u, \varphi, \nabla \varphi ,\mathcal{F}) 
= d(\varphi) Y \left(\nabla u\right)^2 + 
g_c \frac{\gamma}{2} \left(\nabla \varphi\right)^2  
+ \mathcal{K} (\varphi, \mathcal{F}) ,
\end{equation}
where $Y$ represents the Young's modulus, $g_c$, denotes the fracture energy release rate, and $\gamma > 0$ is the phase-field layer width parameter. The degradation function is defined as $d(\varphi) = (1-\varphi)^2$, affecting both elastic response and the electrical conductivity. Additionally, $\mathcal{K} (\varphi, \mathcal{F})$ represents the coupling of damage and fatigue describing how damage evolves due to fatigue over time.

The original model described in \cite{boldrini2016non} is time-dependent. However, our focus is on the long-term behavior of the material, assuming it reaches equilibrium between consecutive time steps. This assumption allows us to simplify the governing equations for $u$ and $\varphi$ into a quasi-static form. However, the evolution of $\mathcal{F}$ that represents long-term aging is still a time-dependent ordinary differential equation (ODE).

\begin{equation}
\label{eq:mechanical}
\nabla \cdot \,\left( (1-\varphi)^2 Y \nabla u \right) 
-  \gamma g_c \, \nabla \cdot \, ( \nabla \varphi \otimes \nabla \varphi  )+  f = 0,
\end{equation}

\begin{equation}
\label{eq:damage}
\gamma g_c \Delta \varphi  
+  (1- \varphi) (\nabla u)^T Y (\nabla u)
- \frac{1}{\gamma}  [ g_c \mathcal{H}' (\varphi) + \mathcal{F} \mathcal{H}_f' (\varphi) ] = 0,
\end{equation}

\begin{equation}
\label{eq:fatigue}
\dot{\mathcal{F}} = -  \frac{ \hat{F}}{\gamma}   \mathcal{H}_f (\varphi),
\end{equation}
The equations (\ref{eq:mechanical}), (\ref{eq:damage}), and 
(\ref{eq:fatigue}) are subjected to appropriate boundary conditions detailed in a later section. The potentials $\mathcal{H}(\varphi)$ and $\mathcal{H}_f(\varphi)$ describe the damage evolution from $0$ to $1$ as fatigue changes from zero to $g_c$. $\mathcal{H}'(\varphi)$ and $\mathcal{H}_f'(\varphi)$ are the derivatives of these potentials with respect to $\varphi$. Appropriate choices of these potentials to ensure a continuous and monotonically increasing transition are:

\noindent\begin{minipage}{.55\linewidth}
	\begin{equation}
	\label{eq:potentials}
	\mathcal{H} (\varphi) =
	\begin{cases}
	0.5 \varphi^2   & \mbox{for} \; 0 \leq \varphi \leq 1 ,
	\vspace{0.1cm} 
	\\
	0.5 + \delta (\varphi -1)  &  \mbox{for} \; \varphi > 1 ,
	\vspace{0.1cm}
	\\
	- \delta \varphi &   \mbox{for} \; \varphi < 0 .
	\end{cases}
	\end{equation} 	
\end{minipage}%
\begin{minipage}{.45\linewidth}
	\begin{equation}
	\label{eq:potentials2}
	\mathcal{H}_f (\varphi) =
	\begin{cases}
	- \varphi & \mbox{for} \; 0 \leq \varphi \leq 1,
	\vspace{0.1cm} 
	\\
	-1 & \mbox{for} \; \varphi > 1,
	\vspace{0.1cm}
	\\
	\hspace{0.25cm} 0 & \mbox{for} \;  \varphi < 0.
	\end{cases}
	\end{equation}
\end{minipage}

We describe the evolution of $\mathcal{F}$ through $\hat{F}$, which represents the formation and growth of micro-cracks under cyclic loading conditions and the influence of temperature. The $\hat{F}$ has a linear relationship with the level of stress associated with the virgin material:

\begin{equation} 
\label{eq:fhat}
\hat{F} =  \rho a \left(\frac{\theta_c}{\theta_0}\right) (1-\varphi) \left|Y \nabla u\right|,
\end{equation}
where the parameter $a$ represents the aging rate, modulated by the ratio of the conductor temperature $\theta_c$ to a reference temperature $\theta_0$ and $\rho$ denotes the material density.

The mechanical model allows damage healing when tensile stress decreases. To simulate an irreversible damage process and avoid healing mechanisms, we adopt an approach similar to that in \cite{miehe2010phase}. We define a variable, $\mathbb{H}$, which represents the local maximum strain energy history:

\begin{equation}
\label{eq:history}
\mathbb{H}(x,t) = \max((\nabla u(x,t))^T Y (\nabla u(x,t)),\mathcal{H}(x,t)).
\end{equation}

We integrate the local maximum strain energy history variable $\mathbb{H}$,  into the damage equation and obtain a modified equation for damage as follows:

\begin{equation}
\label{eq:damage_H}
\gamma g_c \Delta \varphi  
+  (1- \varphi) \mathbb{H}
- \frac{1}{\gamma}  [ g_c \mathcal{H}' (\varphi) + \mathcal{F} \mathcal{H}_f' (\varphi) ] = 0.
\end{equation}

\subsection{Thermal Model}

The original model \cite{boldrini2016non} relates fatigue with the increase in temperature due to repetitive and fast loading; however, here we focus on the long-term damage rather than the short-term increase in temperature. Although the short-term effects are crucial, the long-term reliability of the transmission lines is mostly affected by static loads, environmental conditions, and gradual material degradation. Therefore, we assume a quasi-static regime to simplify our analysis by focusing on steady-state conditions that significantly impact long-term performance. 
Consequently, we adopt the steady heat equation:
\begin{equation}
\label{eq:temperature}
\nabla \cdot (\kappa \nabla \theta_c) + q = 0.
\end{equation}

Equation \ref{eq:temperature} describes the heat balance within a system, where  $\kappa$ represents the thermal conductivity, and  $q$ denotes the net heat exchange. Joule heating $q_j$ due to the current passing through the conductor acts as a heat source while convective cooling $q_c$ due to the wind serves as a heat sink such that the overall heat exchange can be represented as:
\begin{equation}
q = q_j - q_c.
\end{equation}

We consider a convective heat transfer using the following relation.
\begin{equation}
q_c = h (\theta_c - \theta_a),
\end{equation}
where $h$ is the convective heat transfer coefficient. We consider forced convection due to cross flow over the cylinder using the relation given by \cite{cengel2011heat}. The convective heat transfer coefficient is related to the Nusselt number $Nu_D$ and Prandtl number $Pr$ by the following relation: 
 \begin{equation}
Nu_D = \frac{hD}{k} = C Re_D^m Pr^{\frac{1}{3}}.
\end{equation}

We first determine the Reynolds number $Re_D$ using a relation given by 
\begin{equation}
Re_D = \frac{vD}{\nu},
\end{equation}
where $\nu$ is kinematic viscosity of air, $v$ is the velocity of air, and $D$ is the diameter of conductor. 

The values of $C$ and $m$ are experimentally determined which depend upon the value of $Re_D$. Table \ref{tab:experimental values} shows the value of $C$ and $m$ for different ranges of $Re_D$.
\begin{table}[H]
\caption{Values of \(C\) and \(m\) for different \(Re_D\) ranges}
\centering
\begin{tabular}{lll}
\hline
\(Re_D\) Range & \(C\) & \(m\) \\ \hline
0.4 -- 4 & 0.989 & 0.330 \\ 
4 -- 40 & 0.911 & 0.385 \\ 
40 -- 4000 & 0.683 & 0.466 \\ 
4000 -- 40,000 & 0.193 & 0.618 \\ 
40,000 -- 400,000 & 0.027 & 0.805 \\ \hline
\end{tabular}
\label{tab:experimental values}
\end{table}

\subsection{Electrical Model}

In the transmission line design, various factors such as frequency, inductance, reactance, and electromagnetic interactions with the environment and nearby conductors are considered. However, we are not considering the overall design of the transmission lines, rather we focus on understanding how current-driven Joule heating affects the reliability of cables with existing damage. We therefore treat the current $I_b$ as input to our multi-physics system. 
Additionally, we do not consider the transient effects of the AC currents as over time they get averaged out. Therefore, we consider the electric current as a DC-equivalent mean current that remains constant over consecutive time steps. This simplifies our approach for efficient measures of the heat source term. Further to reduce the complexity, we parameterize the value of $I_b$ corresponding to the allowable ampacity value of All Aluminum conductors of diameter around 40 mm. The current data can be obtained through sophisticated methods however, we are not considering a specific transmission line rather we focus on the voltage drop due to temperature and damage-induced resistivity. 

From the above consideration, we focus on solving the conservation of current through the following set of equations:

\begin{equation}
\nabla \cdot J = 0,
\end{equation}

\begin{equation}
J = \sigma_{E} E,
\end{equation}

\begin{equation}
E = -\nabla V,
\end{equation}
where $J$ denotes the electric current density per cross-section area, $E$ is the electric field generated by the voltage $V$, and $\sigma_E$ represents the electric conductivity at the operating temperature, which depends upon the degradation function, $d(\varphi)$ and non-degraded conductivity of damage $\sigma_{E,T}$:

\begin{equation}
\sigma_E = (1 - \varphi)^2 \sigma_{E,T}.
\end{equation}

The non-degraded conductivity $\sigma_{E,T}$ at the operating temperature, is related to the conductivity $\sigma_{E,0}$ at a reference temperature and can be obtained by:

\begin{equation}
\sigma_{E,T} = \frac{\sigma_{E,0}}{1 + \alpha(\theta - \theta_0)},
\end{equation}
where $\alpha$ is the coefficient of resistivity of the conductor. Combining the above equations, we obtain a partial differential equation for the voltage field:

\begin{equation}
\nabla \cdot (-\sigma_{E} \nabla V) = 0,
\label{eq:voltage}
\end{equation}
which can be solved by prescribing either $V$ or $J$ at the boundaries. The voltage field is influenced by damage resulting in a different $\Delta V$ across the conductor.

The Joule heating source $q_j$ can be defined as:

\begin{equation}
q_j = J \cdot E.
\end{equation}

Overall, this model links damage to voltage drop, increasing the power loss due to increased resistance thus aggravating the thermal load due to Joule heating.

\subsection{Sag Consideration}

The temperature of the cable influences the horizontal tension acting on it. This section focuses on obtaining the horizontal load acting on a cable supported by two towers, forming a catenary curve. While several studies have been done to understand the mechanical behavior of such cables \cite{karoumi1999some,stengel2014finite}, we simplify our approach using a one-dimensional damage phase-field model. We assume that the length of the cable $L$ is approximately equal to the span distance $S$. This assumption simplifies the model but does not ignore the presence of sag $D$; rather, it highlights that the primary driver of mechanical damage and fatigue is the horizontal tension $H$.

The sag in the cable is directly related to the temperature variation affecting the horizontal tension. Higher temperatures increase the length of the cable due to thermal elongation, which increases the sag and reduces tension. Conversely, lower temperatures contract the cable reducing the length, ultimately reducing the sag and increasing the tension. Therefore, the horizontal tension at either end of the cable results from an initial pre-tension combined with tension adjustments due to temperature variations. This model determines the appropriate mechanical loading conditions in terms of horizontal tension based on these parameters.

We follow the methodology presented in \cite{grigsby2006electric} to calculate the horizontal tension in the cable. We consider $W_b$ as the weight per unit length and $H_0$ as the initial pre-tension. The initial pre-tension is typically set at about $20\%$ of the material's ultimate strength, to compute the initial sag $S_0$. The initial sag $S_0$ is calculated using:

\begin{equation}
S_0 = \frac{W_b S_L^2}{8 H_0},
\end{equation}

We simplified our model assuming the length of the cable $L$ is equal to the span length $S_L$, despite the actual presence of sag $S$. Although the simplification, the theoretical length is necessary to accommodate the cable with sag over the span:

\begin{equation}
L_0 = S_L + \frac{8 S_0^2}{3 S_L}.
\end{equation}

The length of the cable changes with the change in the temperature according to the classical formula:

\begin{equation}
L = L_0 (1 + \alpha_L \Delta \theta),
\end{equation}

where $\alpha_L$ represents the coefficient of thermal expansion. The resulting change in length alters the sag:

\begin{equation}
S = \sqrt{\frac{3 S_L (L - S_L)}{8}},
\end{equation}

which then leads to a new calculation for the horizontal tension:

\begin{equation}
H = \frac{W S_L^2}{8S},
\end{equation}

In this model, $W$ denotes the total weight, that could account for the weight of additional factors such as ice and wind. For this study, we consider only the wind component and calculate $W$ as:

\begin{equation}
W = \sqrt{W_b^2 + W_w^2}.
\end{equation}

\begin{figure}[H]
    \centering
    \includegraphics[width=0.45\textwidth]{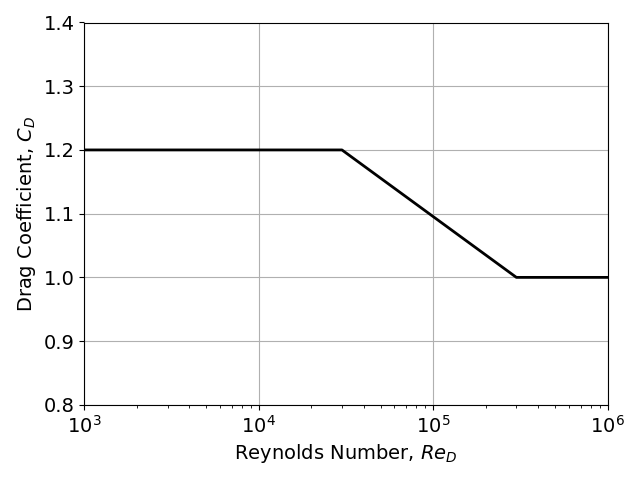}
    \caption{Drag coefficient $C_D$ for the transmission line.}
    \label{fig:CD}
\end{figure}

We follow a similar approach as in \cite{reinoso2020wind,holmes2007wind} to calculate the wind component. The wind-induced component $W_w$ is obtained from the wind pressure $P_w$ and related to the wind velocity $v$:

\begin{equation}
P_w = \frac{1}{2}\rho_{air} v^2,
\end{equation}

\begin{equation}
W_w =  P_w C_D D \sin^2(\theta_w) \alpha
,
\end{equation}

where $\rho_{air}$ is the density of air, $D$ is the diameter exposed to wind, $\theta_w$ is the angle between the transmission line and wind flow, $\alpha$ is the span factor, and $C_D$ is the drag Coefficient. The drag coefficient depends upon Reynold's number as shown in Figure \ref{fig:CD}.

\subsection{Multi-Physics Framework}
\begin{figure}[H]
    \centering
    \includegraphics[width=0.45\textwidth]{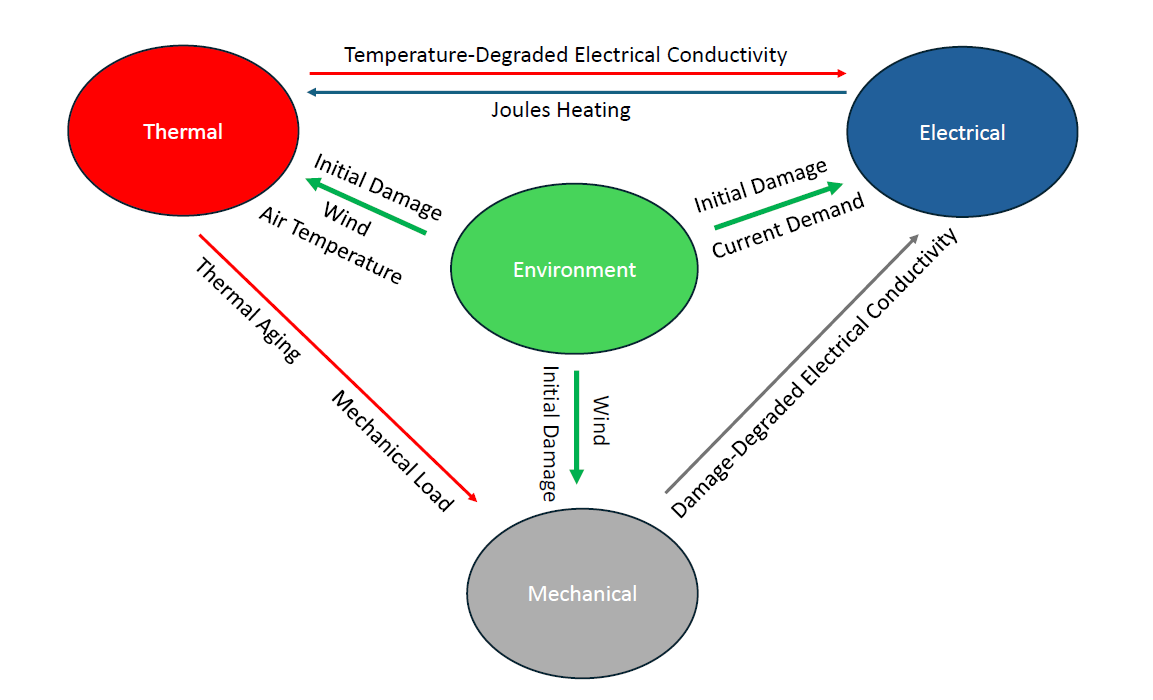}
    \caption{\textit{Schematic diagram illustrating the relationship between the thermal, electrical, and mechanical models, in addition to an environmental module responsible for initial conditions and input parameters.}}
    \label{fig:relations}
\end{figure}

In Figure~\ref{fig:relations}, we present a schematic representation outlining our approach to modeling transmission line failures. The model consists of the main governing equations, including Eqs (\ref{eq:mechanical}), (\ref{eq:damage_H}), (\ref{eq:fatigue}), (\ref{eq:temperature}), and (\ref{eq:voltage}), which describe the interactions between the thermal, electrical, and mechanical aspects of the system. The diagram offers a detailed view of how the thermal, electrical, and mechanical components are interconnected within the system. Additionally, we introduce an Environmental module, which is responsible for setting the inputs and initial conditions for all other modules based on the specific scenario under consideration and its relevant parameters.

\section{Deterministic Solution}
In this section, we explore the finite element discretization within our multi-physics framework and outline the deterministic solution process. This deterministic approach has two primary objectives: it acts as a black box for non-intrusive stochastic analysis for uncertainty quantification, and it provides the interpretable basis for evaluating multi-dimensional uncertainty propagation.
\subsection{Finite-Element Discretization}
In a one-dimensional setting, variables are treated as scalar fields, simplifying the representation of differential operators. We use the standard linear finite-element method for spatial discretization. By multiplying the governing equations—Eqs.(\ref{eq:mechanical}), (\ref{eq:damage_H}), (\ref{eq:fatigue}), (\ref{eq:temperature}), and (\ref{eq:voltage})—by a test function 
$w$, integrating by parts, and considering the volume differential expressed as the cross-sectional area 
$A(x)$, we derive their corresponding weak forms as follows
\begin{equation}
\int_0^L -(1 - \varphi)^2 Y A(x) \frac{du}{dx}\frac{dw}{dx} dx + \int_0^L \gamma g_c A(x) \left(\frac{d\varphi}{dx}\right)^2 \frac{dw}{dx} dx + \int_0^L f A(x) w dx = 0,
\end{equation}

\begin{equation}
\begin{split}
\int_0^L - \gamma g_c A(x) \frac{d\varphi}{dx}\frac{dw}{dx} dx + \int_0^L A(x) \mathbb{H} w dx - \int_0^L A \mathbb{H} \varphi w dx \\ - \int_0^L \frac{g_c A(x)}{\gamma} \varphi w dx + \int_0^L \frac{A(x)}{\gamma}\mathcal{F}w dx = 0,
\end{split}
\end{equation}

\begin{equation}
\int_0^L \dot{\mathcal{F}} w A(x) dx = \int_0^L \frac{-\rho a (1 - \varphi) Y \vert\frac{du}{dx}\vert (-\varphi)}{\gamma}\frac{\theta_c}{\theta_0} w A(x) dx, 
\end{equation}

\begin{equation}
\begin{split}
\int_0^L - \kappa A(x) \frac{d\theta}{dx}\frac{dw}{dx} dx + \int_0^L \sigma_E A(x)\left(\frac{dV}{dx}\right)^2 w dx \\- \int_0^L h \theta A_s(x) w dx + \int_0^L h \theta_{a} A_s(x) w dx = 0,
\end{split}
\end{equation}

\begin{equation}
\int_0^L \sigma_E \frac{dV}{dx} \frac{dw}{dx} A(x) w dx = 0,
\end{equation}
where we define $A_s(x)$ as the variable surface area over which convective heat transfer due to wind occurs.

For each element $k$, we use a linear approximation where the field variables are represented as linear combinations of nodal basis functions:

\begin{align}
u^k &= N \hat{u}^k,\\
\varphi^k &= N \hat{\varphi}^k ,\\
\mathcal{F}^k &= N \hat{\mathcal{F}}^k ,\\ 
\theta^k &= N \hat{\theta}^k ,\\
V^k &= N \hat{V}^k.
\end{align}

We calculate the finite-element interpolation of spatial derivatives using linear combinations of the derivatives of shape functions:

\begin{align}
\left(\frac{du}{dx}\right)^k &= B \hat{u}^k,\\
\left(\frac{d\varphi}{dx}\right)^k &= B \hat{\varphi}^k,\\
\left(\frac{d\theta}{dx}\right)^k &= B \hat{\theta}^k,\\
\left(\frac{dV}{dx}\right)^k &= B \hat{V}^k,\\
\end{align}
where we define $N$, $B$, $\hat{u}^k$, $\hat{\varphi}^k$, $\hat{\mathcal{F}}^k$, $\hat{\theta}^k$, $\hat{V}^k$ as
\begin{align}
N &= \begin{bmatrix}
N_1 & N_2
\end{bmatrix},\\
B &= \begin{bmatrix}
N_{1,x} & N_{2,x}
\end{bmatrix},\\
\hat{u}^k &= \begin{bmatrix}
u_1^k & u_2^k
\end{bmatrix},\\
\hat{\varphi}^k &= \begin{bmatrix}
\varphi_1^k & \varphi_2^k
\end{bmatrix},\\
\hat{\mathcal{F}}^k &= \begin{bmatrix}
\mathcal{F}_1^k & \mathcal{F}_2^k
\end{bmatrix},\\
\hat{\theta}^k &= \begin{bmatrix}
\theta_1^k & \theta_2^k
\end{bmatrix},\\
\hat{V}^k &= \begin{bmatrix}
V_1^k & V_2^k
\end{bmatrix},
\end{align}
where $N_1$ and $N_2$ are linear interpolation functions.

We substitute the earlier approximations into the weak form and use a forward Euler method to evolve $\mathcal{F}$, resulting in the following discretization for each $k$th element: 

\begin{align}
K_u \hat{u}^k &= w_u + M \hat{f}^k,\\
K_{\varphi} \hat{\varphi}^k &= w_\varphi,\\
M \hat{\mathcal{F}}^{{n+1}^k} &= M \hat{\mathcal{F}^{n}}^k + \Delta t w_{\mathcal{F}},\\
K_\theta \hat{\theta}^k &= w_\theta,\\
K_V \hat{V}^k &= 0.
\end{align}
Here, the superscripts $n$ and $n+1$ denote the current and subsequent time steps, respectively. The discrete forms are defined using these operator conventions:

\begin{align}
K_u &= \int_k (1 - N\hat{\varphi}^k)^2 Y A(x) B^T B \, dx,\\
w_u &= \int_k \gamma g_c A(x) (B \hat{\varphi}^k)^2 B \, dx,\\
M &= \int_k A(x) N^T N \, dx,\\
K_\varphi &= \int_k \gamma g_c A(x) B^T B \, dx + \int_k \mathbb{H} A(x) N^T N \, dx + \int_k \frac{g_c A(x)}{\gamma} N^T N \, dx,\\
w_\varphi &= \int_k \mathbb{H} A(x) N \, dx + \int_k \frac{A(x)}{\gamma} N^T \hat{\mathcal{F}^n}^k N \, dx,\\
K_\theta &= \int_k \kappa A(x) B^T B \, dx + \int_k h A_s(x) N^T N \, dx,\\
w_\theta &= \int_k \sigma_E A(x) \left(B \hat{V}^k\right)^2 N \, dx  + \int_k h A_s(x) \theta_{a} N \, dx,\\
K_v &= \int_k (1 - N\hat{\varphi}^k)^2 \sigma_{E} A(x) B^T B \, dx.
\end{align}

We derive the global forms of these matrices and vectors using standard finite-element assembly procedures. 
\begin{align}
K_u \hat{u} &= w_u + M \hat{f}, \label{eq:u}\\
K_{\varphi} \hat{\varphi} &= w_\varphi, \label{eq:phi}\\
M \hat{\mathcal{F}}^{n+1} &= M \hat{\mathcal{F}}^n  + \Delta t w_{\mathcal{F}}, \label{eq:F}\\
K_\theta \hat{\theta}  &= w_\theta, \label{eq:theta}\\
K_V \hat{V} &= 0 \label{eq:V}.
\end{align}

The formulation supports a staggered solution scheme at each time step, implemented according to the following algorithm:
\begin{algorithm}[H]
\caption{Solution of Thermo-Electro-Mechanical Model.}
\label{algo:deterministic}
\begin{algorithmic}[1]
\State Choose initial pre-tension.
\For{Each time-step}
\State Compute the current tensile load.
\State Solve for displacements.
\State Update strain energy history.
\State Solve damage field.
\State Update fatigue.
\State Solve the temperature field.
\State Solve voltage field.
\EndFor
\end{algorithmic}
\end{algorithm}

\subsection{Numerical Results}
We consider an all aluminum-type conductors are subjected to cyclic loading conditions under wind, temperature, and current. The details on wind and temperature loading conditions are provided in Section~\ref{sec:problem statement}. The current loading is parameterized using the following relation:
\begin{equation}
    I(t) = -I_b - I_a (\sin 4 \pi t),
\end{equation}
where $I_b$ is the base current, set at 1500 A, and $I_a$ is the amplitude, set at 100 A, corresponding to the allowable ampacity for a 40 mm aluminum conductor.

We represent the horizontal tension by setting $u=0$ at $x=0$ and $H$ at $x=L$. For damage, we specify $\frac{d \varphi}{dx} = 0$ at both boundaries. Similarly, for the current conservation equation, we apply boundary conditions similar to the mechanical case: setting $V = 0$ at $x=0$ and applying a current density $J$ at $x=L$.

\begin{figure}[H]
    \centering
    \includegraphics[width=0.45\textwidth]{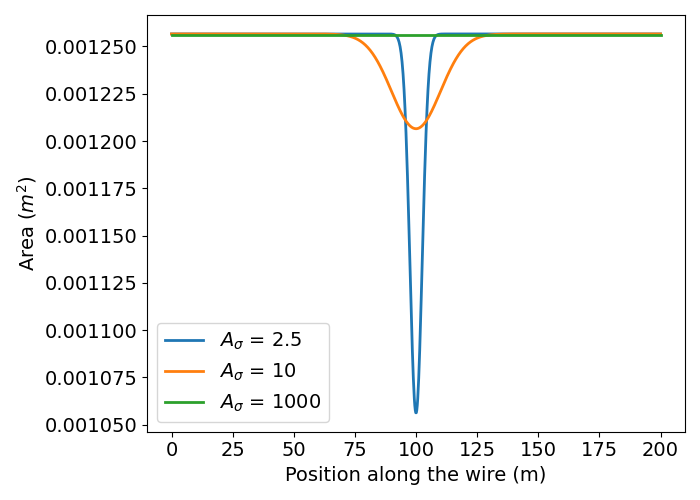}
    \caption{Cross-section area variation as a function of $A_\sigma$}
    \label{fig:Area}
\end{figure}

Additionally, we represent the initial damage as a variable cross-section area across the line. In practice, all materials feature inherent imperfections that accumulate over time until a critical point where significant damage occurs. To simulate this, we assume a reduced cross-section area at the center of the cable, representing the cumulative effects of multiple defects and damage that might lead to a fracture. We also consider a base scenario where the materials are assumed to have insignificant imperfections, serving as a comparison point with the significant initial damage cases. We define the cross-section area using the following relation:
\begin{equation}
A(x) = A_0 \left(1 - \frac{1}{A_\sigma \sqrt{2 \pi}} \exp\left( \frac{- (x - L/2)^2 }{2 A_\sigma^2}\right)\right),
\end{equation} 
where $A_0 $ denotes the undamaged cross-section area, while $A_\sigma$ indicates the ratio of spread to depth of the variation in the area, representing different levels of damage. Figure~\ref{fig:Area} illustrates various area profiles based on $A_\sigma$.
The material parameter values of aluminum conductor are provided in Table~\ref{tab:param}. All simulations are conducted to a 6000 run period with the time-step $\delta t$ = 0.01, corresponding to a life-cycle of 60 years. The aluminum conductors start to anneal at temperatures above 366 K \cite{vasquez2017end,hathout2018impact}, and rupture if the temperature exceeds 373 K \cite{cimini2013temperature}. Therefore, we set the maximum temperature limit to $\theta_{\text{lim}} = 373$ K (100°C). Beyond this limit, all simulations are stopped.

\begin{table}[H]
\caption{Air, Material, and Geometry parameters.}
    \centering
    \begin{tabular}{lll} 
        \hline
        \textbf{Parameter} & \textbf{Value} & \textbf{Unit} \\ 
        \hline
        Length of Cable $L$ & 200 & m \\
        Number of Elements $N$ & 1000 &  \\
        Diameter $D$ & 0.04 & m \\
        \hline
        Young modulus $Y$ & 69 & GPa \\
        Damage layer width $\gamma$ & 0.02 & m \\
        Fracture energy $g_c$ & 10 & kN/m \\
        Density $\rho$ & 2700 & kg/m$^3$ \\
        Aging coefficient $a$ & $1 \times 10^{-10}$ & m$^5$/(y kg) \\
        Thermal conductivity $\kappa$ & 237 & W/(m K) \\
        Electrical conductivity $\sigma_{E,0}$ & $3.77 \times 10^7$ & S/m \\
        Temperature coefficient $\alpha$ & $3.9 \times 10^{-3}$ & K$^{-1}$ \\
        \hline
        Density of air $\rho_{air}$ & 1.225 & kg/m$^3$ \\
        Kinematic viscosity of air $\nu$ & $15 \times 10^{-6}$  &  m$^2$/s\\
        Thermal Conductivity of air $\kappa_{air}$ & 0.0295  &  W/(m K)\\
        Prandtl Pr & 0.71 &  \\
        \hline
    \end{tabular}
    \label{tab:param}
\end{table}

\begin{figure}[H]
	\centering
	\subfloat[Damage.]{\includegraphics[width=0.25\textwidth]{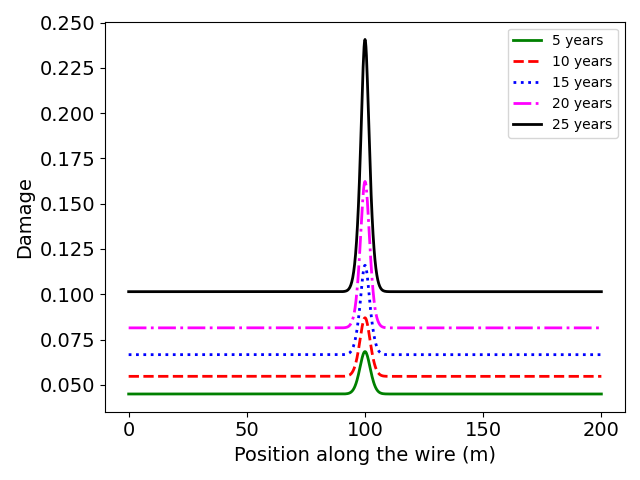}}
	\subfloat[Fatigue.]{\includegraphics[width=0.25\textwidth]{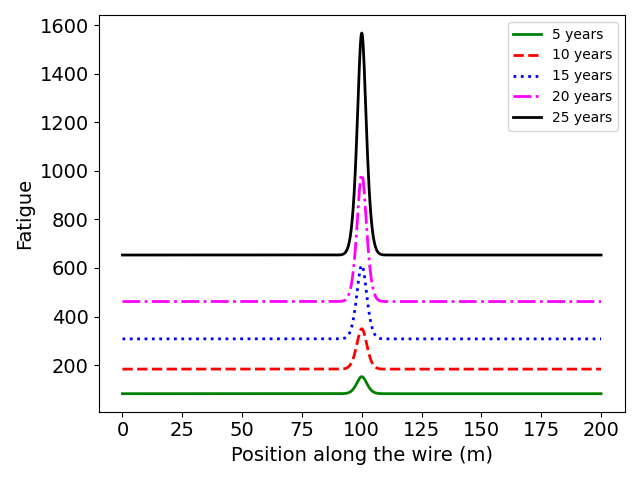}}
	\subfloat[Temperature.]{\includegraphics[width=0.25\textwidth]{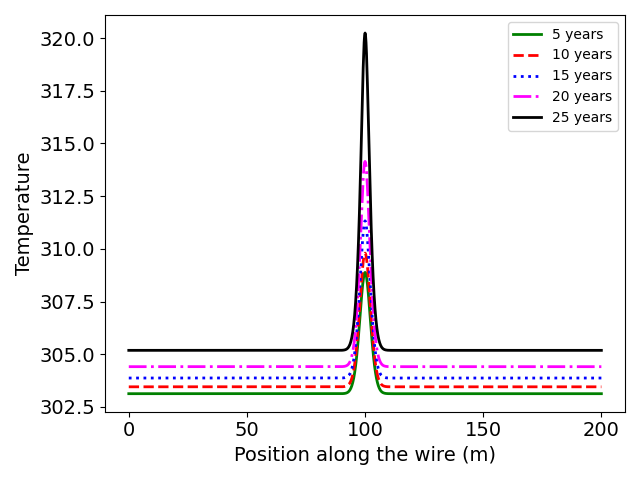}}
	\subfloat[Voltage drop.]{\includegraphics[width=0.25\textwidth]{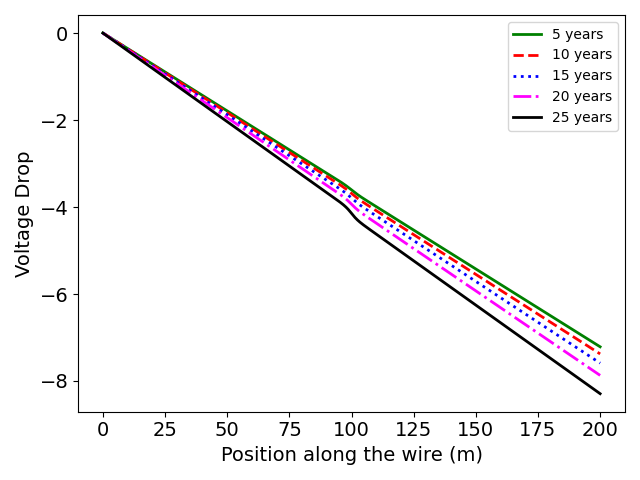}}
	\caption{Evolution of field variables.}
	\label{fig:field_1_1}
\end{figure}

We begin our analysis by examining the evolution of field quantities in the Texas scenario, using the reference parameters listed in Table~\ref{tab:param}. Every five years, we plot the evolution of damage, fatigue, temperature, and voltage fields along the length of the transmission line as shown in Fig~\ref{fig:field_1_1}. We note that damage typically starts and accumulates in regions with smaller cross-sectional areas, leading to increased temperatures and notable disturbances in the voltage fields. Additionally, as damage and temperature increase, the voltage drop along the line increases over time due to increased electrical resistance in the conductor.

\begin{figure}[H]
	\centering
	\subfloat[Damage.]{\includegraphics[width=0.25\textwidth]{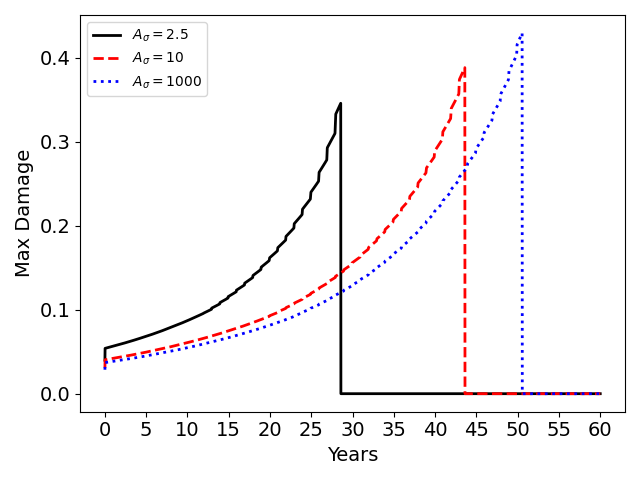}}
	\subfloat[Fatigue.]{\includegraphics[width=0.25\textwidth]{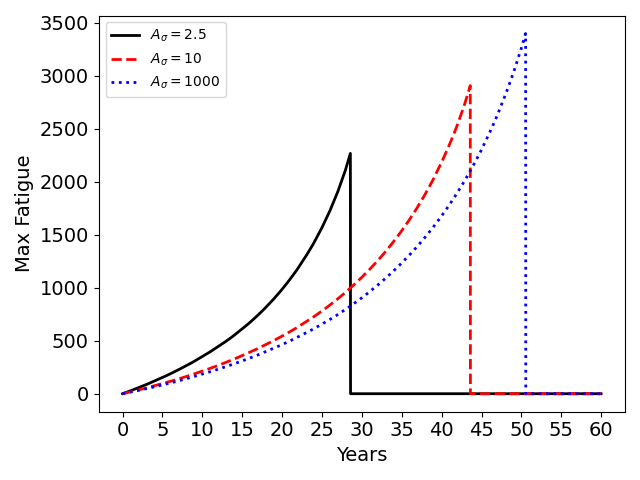}}
	\subfloat[Temperature.]{\includegraphics[width=0.25\textwidth]{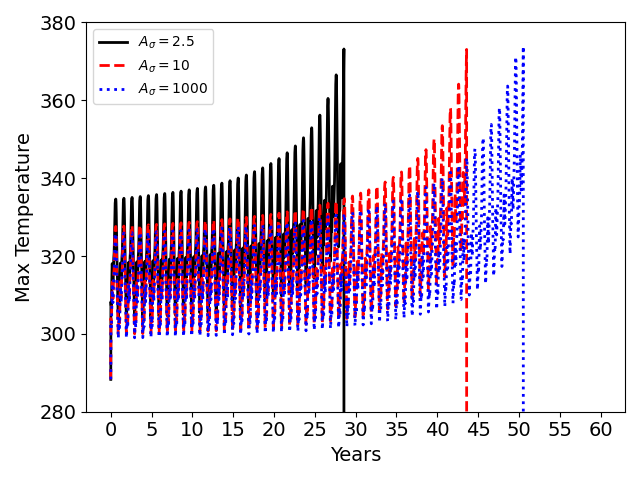}}
	\subfloat[Voltage drop.]{\includegraphics[width=0.25\textwidth]{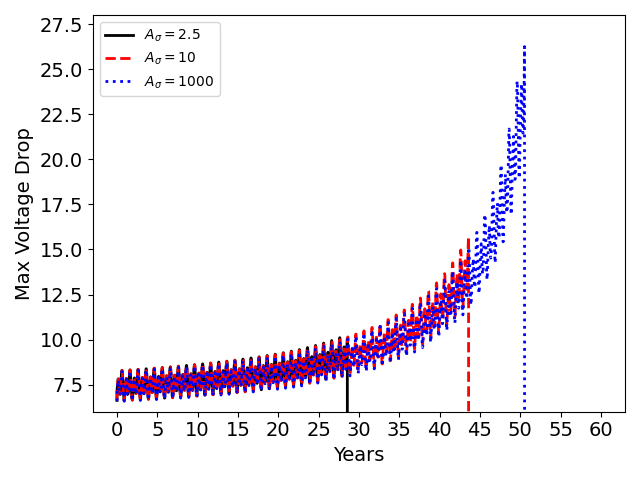}}
	\caption{Effect of initial damage on maximum field values over time.}
	\label{fig:ts_1_1}
\end{figure}
Next, we study the time evolution of maximum values under the effect of initial damage $A_{\sigma}$, by tracking the central node in cases of damage, fatigue, and temperature, as well as the maximum absolute voltage drop from the end node in the voltage field in Figure \ref{fig:ts_1_1}. 
\begin{figure}[H]
    \centering
    \subfloat[Texas]
    {\includegraphics[width=0.45\textwidth]{ts_temperature_1_123.png}}
    \subfloat[California]{\includegraphics[width=0.45\textwidth]{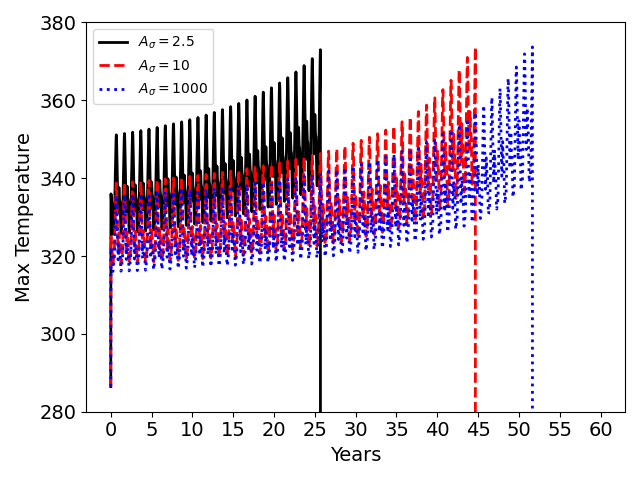}}
    \hfill 
    \subfloat[Michigan]
    {\includegraphics[width=0.45\textwidth]{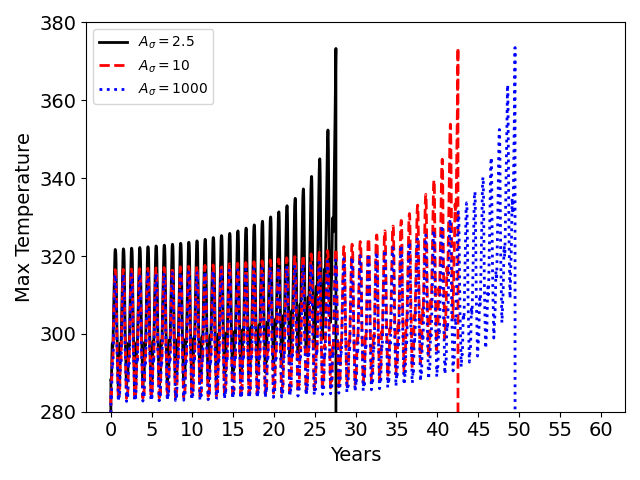}}
    \subfloat[Florida]
    {\includegraphics[width=0.45\textwidth]{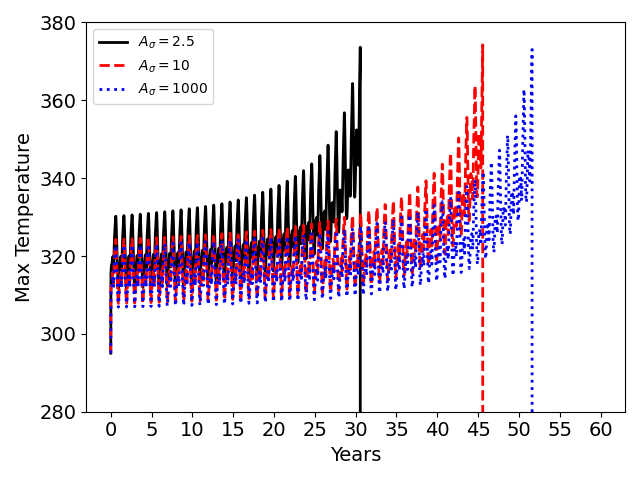}}
    \caption{Failure of Transmission line for different values of initial damage}
    \label{fig:deterministic}
\end{figure}

Moreover, we compare the life span of transmission lines for each specific scenario varying the cross-sectional area of damage $A_{\sigma}$ under respective loading conditions in Figure\ref{fig:deterministic}. In the Texas scenario, we observe the life span of the transmission line reduces from 51 years to 44 years with moderate damage. However, with severe damage, the life span reduces to 28 years indicating how detrimental initial damage could be on the longevity of the material. Similarly, in the California scenario, the lifespan of the transmission line reduces from 52 years to 26 years under severe damage. This reduction is attributed to the state's high temperatures and low wind speeds, which significantly lower the convective cooling effect. Consequently, in the absence of the cooling effect, the accumulation of damage raises the material temperature, leading to overheating and early failure. In the case of the Michigan scenario, the state shows the most fluctuating temperature and wind over a year which progresses the damage raising the temperature of the material and leading to early failure under insignificant damage conditions among the four states. However, due to the high wind speed and low temperature enhancing convective cooling during some seasons, the life span under severe damage is higher than the California scenario. The Florida scenario shows high temperatures over a year similar to Texas and California, however, the high wind speed in the state enhances the convective cooling reducing the early failure of the transmission line in the state. Overall, every scenario shows a significant reduction in the life span of the transmission line under severe damage.

\section{Stochastic Solution}
In this section, we address the challenge of performing uncertainty quantification (UQ), sensitivity analysis (SA), and probability of failure ($p_f$) on the model treated as a black box.  We focus on the conductor temperature as our QoI. The non-intrusive nature of our approach is particularly beneficial, allowing us to apply the same procedures utilized in deterministic solutions to each realization of the stochastic problem without altering the model's governing equations.

We employ the PCM for three main purposes. First, we use PCM to compute the moments of our QoI. Second, we conduct global sensitivity analysis by calculating Sobol Sensitivity Indices ($S_i$). These indices highlight the relative importance of each model parameter by quantifying the contribution to the total variance of the QoI. The computational efficiency of PCM, especially for the calculation of $S_i$, simplifies the sensitivity analysis by post-processing UQ data in a straightforward manner. 
Lastly, we use PCM to facilitate the computation of the probability of failure directly from the first moment of the Bernoulli random variable defined directly from the limit state function.

\subsection{Uncertainty Quantification}
To conduct the uncertainty quantification (UQ) analysis, we employ the PCM, which utilizes polynomial interpolation to approximate solutions within the stochastic space. The method involves mapping the points from the physical space to the stochastic space using the parametric probability density function (PDF) by approximating the solution using orthogonal Lagrange polynomials. Due to their orthogonal properties, these polynomials simplify the computation of expectations and variances by evaluating the QoI at the collocation points significantly reducing computational costs and enhancing convergence rates. 

Following the methodology outlined in \cite{barros2021integrated}, let $(\Omega_s,\mathcal{G}, \mathbb{P})$ represent a complete probability space, where $\Omega_s$ denotes the space of outcomes $\omega$, $\mathcal{G}$ is the $\sigma-$algebra, and $\mathbb{P}$ is a probability measure mapping $\mathcal{G}$ to the interval [0,1]. In our model, we treat the material and loading parameters as random variables with two different sets as, $\xi_m(\omega)$ and $\xi_l(\omega)$. We model our QoI temperature as random variables. For simplicity, we denote the random parameters by $\xi = \xi(\omega)$.

Our quantity of interest, denoted as $Q$, involves calculating the mathematical expectation $\mathbb{E}[Q(x, t; \xi)]$ within a one-dimensional stochastic space as
\begin{equation}
\label{eq:pcm_expectation}
\mathbb{E}\left[ Q(x, t; \xi) \right] = \int_{a}^b Q(x, t; \xi) \rho(\xi) d \xi,
\end{equation}
\noindent where PDF of $\xi$, denoted as $\rho(\xi)$, is evaluated by mapping the physical parametric space to the standard space $[-1, 1]$ and integrating using  Gauss quadrature. This transformation allows the integral to be expressed within this standard interval as follows
\begin{equation}
\label{eq:pcm_expectation2}
\mathbb{E}\left[ Q(x, t; \xi) \right] = \int_{-1}^1 Q(x, t; \xi(\eta)) \rho(\xi(\eta)) J d \xi(\eta),
\end{equation}

\noindent where $J = d\xi / d\eta$ represents the Jacobian of the transformation. To approximate the expectation, we use polynomial interpolation in the stochastic space, expressed as $\hat{Q}(x, t; \xi)$, which approximates the exact solution as

\begin{equation}
\label{eq:pcm_expectation3}
\mathbb{E}\left[ Q(x, t; \xi) \right] \approx \int_{-1}^1 \hat{Q}(x, y,t; \xi(\eta)) \rho(\xi(\eta)) J d \xi(\eta).
\end{equation} 	
We use Lagrange polynomials $L_i(\xi)$ to interpolate the solution in the stochastic space:

\begin{equation}
\label{eq:poly_approximation}
\hat{Q}(x, t; \xi) = \sum_{i = 1}^{I} Q(x, t; \xi_i) L_i(\xi),
\end{equation}

\noindent which satisfies the Kronecker delta property at the interpolation points:

\begin{equation}
\label{eq:delta}
L_i(\xi_j) = \delta_{ij}.
\end{equation}
We substitute the polynomial approximation from Eq.(\ref{eq:poly_approximation}) into Eq.(\ref{eq:pcm_expectation3}), and used the quadrature rule to approximate the integral and compute the expectation as

\begin{equation}
\mathbb{E}\left[ Q(x, t; \xi) \right] \approx \sum_{p = 1}^{P} w_p \rho(\xi(\eta)) J \sum_{i = 1}^{I} Q(x, t; \xi(\eta) ) L_i(\xi(\eta) ), \label{eq:pcm_exp_approx2}
\end{equation}
where $\eta_p$ and $w_p$ are coordinates and weights, respectively, for each integration point $q = 1,\, 2,\, \dots,\, P$. We choose the collocation points to be the same as the integration points in the parametric space by the Kronecker property of the Lagrange polynomials in Eq.(\ref{eq:delta}) and simplify the approximation from Eq.(\ref{eq:pcm_exp_approx2}) into a single summation:

\begin{equation}
\label{eq:pcm_single}
\mathbb{E}\left[ Q(x, t; \xi) \right] = \sum_{p = 1}^{P} w_p \rho(\xi_p(\eta_p))  J Q(x, t; \xi_p(\eta_p)).
\end{equation}
We use a linear affine mapping from the standard domain to the real domain, defined as $\xi_p(\eta_p) = a + \frac{(b-a)}{2}(\eta_p + 1)$. This mapping calculates the Jacobian for one-dimensional integration as $J = (b-a)/2$. It also provides the values of the random variable in the physical space.

Finally, we approximate the integration and rewrite as a summation over the collocation points, assuming a uniform distribution for the parameters throughout the interval $[a, b]$, with $\rho(\xi) = 1/(b-a)$. The expectation becomes 
\begin{equation}
\mathbb{E}\left[ Q(x, y, t; \xi) \right] = \frac{1}{2}\sum_{p = 1}^{P}  w_p Q(x, t; \xi_p).
\end{equation} 

Similar to the MC method, the standard deviation is computed as

\begin{equation}
\sigma \left[ Q(x, t; \xi) \right] =  \sqrt{ \frac{1}{2}\sum_{p = 1}^{P}  w_p \left( Q(x, t; \xi_p) - \mathbb{E}\left[ Q(x, t; \xi) \right] \right)^2}.
\end{equation}

The generalization of the PCM to higher dimensions involves adding additional integrals to Eq.~(\ref{eq:pcm_expectation}) which then reduces to

\begin{align}
\mathbb{E}\left[ Q(x, t; \xi^1,\, \dots ,\, \xi^k) \right] &= \mathbb{E}_{PCM}\left[ Q(x, t; \xi^1,\, \dots ,\, \xi^k) \right] \notag \\ &\approx \sum_{p = 1}^{P}\dots \sum_{l = 1}^{L} w_p \dots w_l \, \rho(\xi_p) \dots \rho(\xi_l) \, J_p \dots J_l\,  Q(x,  t; \xi^1_p,\,\dots,\, \xi^k_l) \label{eq:pcm_multi}
\end{align}

In the expanded version of PCM for higher dimensions, we have $k$ summations, corresponding to each dimension in the random space. Here, $\xi^k_l$ is notated such that the superscript denotes the dimension in the random space, and the subscript indicates the specific collocation point within that dimension. For simplicity, we use $\mathbb{E}\left[ Q(x, t; \xi^1,, \dots ,, \xi^k) \right] = \mathbb{E}\left[ Q\right]$ to represent the expectation. We then formulate the expression for the standard deviation as follows:

\begin{align}
&\sigma\left[ Q(x, t; \xi^1,\, \dots ,\, \xi^k) \right] = \sigma_{PCM}\left[ Q(x, t; \xi^1,\, \dots ,\, \xi^k) \right] \notag \\&\approx \sqrt{\sum_{p = 1}^{P}\dots \sum_{l = 1}^{L} w_p \dots w_l \, \rho(\xi_p) \dots \rho(\xi_l) \, J_p \dots J_l\, \left(Q(x, t; \xi^1_p,\,\dots,\, \xi^k_l) -  \mathbb{E}\left[ Q\right]\right)^2}. \label{eq:pcm_multi_std}
\end{align}

We assume that the random variables are mutually independent and the discretization in the parametric space is isotropic. We further remark the fully tensorial product used in PCM is sufficient for the model in this study, as it involves fewer than 6 dimensions. However, for models with more than 6 dimensions, this approach becomes inefficient due to the exponential increase in the number of simulations required by the tensor product. To reduce the problem of dimensionality in higher-dimensional stochastic spaces, Smolyak sparse grids \cite{smolyak1963quadrature} can be used to reduce the number of realizations while maintaining accuracy. Additional techniques for reducing dimensionality in UQ include Principal Component Analysis (PCA) \cite{abdi2010principal}, low-rank approximations \cite{chevreuil2015least}, and active subspaces methods \cite{constantine2017global}.

\subsection{Global Sensitivity Analysis}
We examine the global sensitivity of input parameters using Sobol indices \cite{sobol1993sensitivity}, which help us to determine the relative importance of each parameter to the variance of our QoI. For details on the derivation, we refer to the work by Saltelli et al. (2010) \cite{saltelli2010variance}. In our global sensitivity analysis, each parameter is denoted as $\xi^j$, where $j = 1, 2, ..., k$, and $k$ represents the total number of dimensions in the parametric space. We evaluate the influence of the parameter $\xi^j$ on the variance $V$ of the QoI as follows:
\begin{equation}
\label{eq:var}
V_{\xi^j}\left(\mathbb{E}_{\mathbf{\xi}^{\sim j}}(Q | \xi^j)\right)
\end{equation}
\noindent where $\mathbf{\xi}^{\sim j}$ represents all possible values of the random parameters except for $\xi^j$, which is held fixed. Equation~(\ref{eq:var}) describes the process of computing the expected value of $Q$ with $\xi^j$ fixed, followed by calculating the variance across all possible values of $\xi^j$. Based on the Law of Total Variance, we have:

\begin{equation}
\label{eq:law}
V_{\xi^j}\left(\mathbb{E}_{\mathbf{\xi}^{\sim j}}(Q | \xi^j)\right) + \mathbb{E}_{\mathbf{\xi}^{j}}\left(V_{\mathbf{\xi}^{\sim j}}(Q | \xi^j)\right) = V(Q)
\end{equation}

The second term on the left-hand side is residual and $V(Q)$ is the total variance. We normalize Eq.~(\ref{eq:law}) to obtain the first-order sensitivity index. This index quantifies the impact of the random variable $\xi^j$ on the total variance, and is calculated as follows:

\begin{equation}
\label{eq:si}
S_i = \frac{V_{\xi^j}\left(\mathbb{E}_{\mathbf{\xi}^{\sim j}}(U | \xi^j)\right)}{V(U)}
\end{equation}

The sensitivity indices, $S_i$, quantify the first-order effects of the variable $\xi^j$ on the variance, excluding interactions between $\xi^j$ and other parameters. Following normalization, the sum of all $S_i$ is less than 1 with the remaining portion representing higher-order interactions among the parameters, which are not considered in this paper but could be analyzed similarly through post-processing of the PCM.

The computation of $S_i$ could be challenging when using methods like Monte Carlo (MC) for uncertainty quantification (UQ) due to their computational intensity. However, in this case, the PCM serves as an effective building block for fast and cost-efficient computations of global sensitivity.

\subsection{Probability of Failure}
In the final phase of our stochastic analysis on transmission line failure, we focus on calculating the probability of failure ($p_f$) over time. Traditionally, MC methods are employed in reliability analysis to compute the probability of failure \cite{machado2015reliability}. Stochastic collocation methods are used for similar analysis, however, the method primarily yields the moments of a limit state function  $g(R, S)$. To compute 
$p_f$, these moments must be transformed into a probability density function (PDF), which can be achieved using various techniques such as the method of moments \cite{low2013new,dang2019novel}, Polynomial chaos \cite{lasota2015polynomial,garcia2021polinomial}, Gaussian transformations \cite{he2014sparse}, or entropy optimization methods \cite{winterstein2013extremes}. However, these processes add a layer of complexity.  

In this study, we present an alternative approach that utilizes the efficiency of the PCM to compute the probability of failure ($p_f$). Rather than calculating the moments of the limit state function $g$, obtaining an approximate PDF, and then determining $P(g<0)$, we simplify the process by transforming $g$ into a Bernoulli random variable $h_B$ with coefficient $p_h$.

The transformation of $g$ to $h_B$ is defined as follows:
\begin{equation}
\label{eq:h}
h_B = \begin{cases}
0, & \text{if } g \geq 0,\\
1, & \text{otherwise}.
\end{cases}
\end{equation}

In practice, each realization of the PCM generates a time-series vector $h_B$. It remains zero until the point where $\theta_{\text{max}}$ exceeds the threshold $\theta_{\text{lim}}$, at which $h_B$ becomes one and remains at this value until the final time-step. This represents a step function at the point of failure for each realization.

At a fixed time-step, considering the expectation of $h_B$, due to the smoothness of the QoI, $h_B$ can be assumed as a real value between 0 and 1 which reflects the probability of the maximum temperature exceeding the limit up to that time.

Using PCM to compute the expectation of $h_B$, which behaves as a Bernoulli random variable, provides $p_f$, the probability of failure. This expected value of $h_B$ directly represents the Bernoulli parameter $p_h$. Thus, a single integration using PCM effectively and accurately measures $p_f$ for each time step.

\subsection{Numerical Results}
In deterministic simulations, we have a solid understanding of how varying initial damage under cyclic loading conditions influences the lifespan of the transmission line. Now, we shift our focus to the impact of parametric uncertainty on maximum temperature. For this analysis, we model the parametric uncertainty with a uniform distribution, where each parameter varies by $10\%$ around its mean value, corresponding to the parameters used in the deterministic model.

\subsubsection{Texas Scenario}
We initially consider the Texas Scenario as our baseline for conducting preliminary analyses on the uncertainties associated with parametric material properties and input loading conditions. Our primary goal is to quantify uncertainty and perform sensitivity analysis on factors influencing the maximum temperature. We focus on the uncertainties in material parameters, represented by the set $\xi_m(\omega) = \{A_\sigma(\omega), \gamma(\omega), g_c(\omega), a(\omega)\}$. These parameters are selected due to their inherent measurement inaccuracies or their assumptions within mathematical modeling. All other material parameters are considered deterministic. Through UQ and SA, we aim to identify the two parameters within $\xi_m(\omega)$ that exert the most significant impact on the variance of $\theta_{\text{max}}$.
Following the analysis of material parameters, we then examine the uncertainty in loading conditions, which are represented within a separate stochastic space, denoted as $\xi_l(\omega) = {\theta_b(\omega), w_b(\omega), I_b(\omega), I_A(\omega)}$.
From this set, we identify three parameters that are most influential according to global SA results.

In all our simulations, we use 5 PCM points per dimension. It is important to note that when analyzing time-series data for maximum temperature across the entire stochastic space, we must adjust for the varying times of failure observed in each realization. To manage this, we truncate the time series at the earliest time of failure across all realizations.

\begin{figure}[H]
	\centering
	\subfloat[Expectation.]{\includegraphics[width=0.45\textwidth]{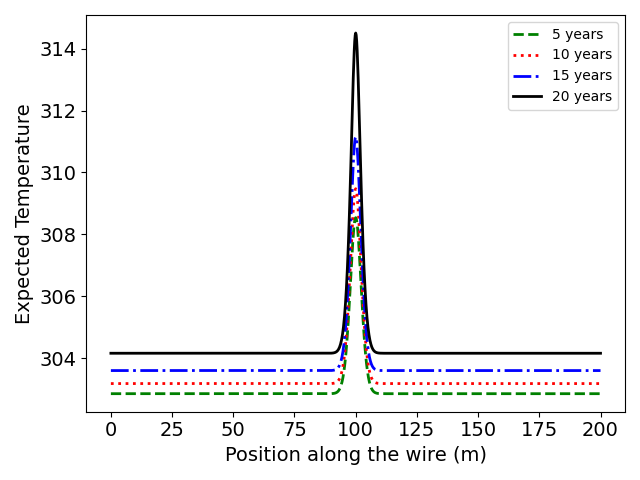}}
	\subfloat[Standard deviation.]{\includegraphics[width=0.45\textwidth]{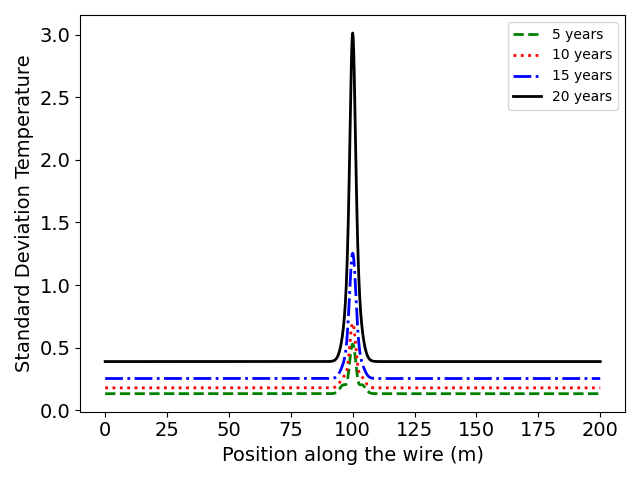}}
	\caption{Expectation and standard deviation temperature fields for material parameter uncertainty set $\xi_m(\omega)$ in the Texas scenario.}
	\label{fig:uq_1_1}
\end{figure}

\begin{figure}[H]
	\centering
	\subfloat[Expectation.]{\includegraphics[width=0.45\textwidth]{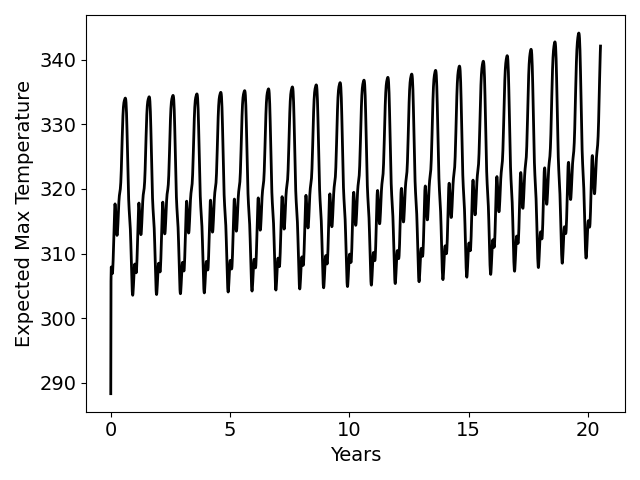}}
	\subfloat[Standard deviation.]{\includegraphics[width=0.45\textwidth]{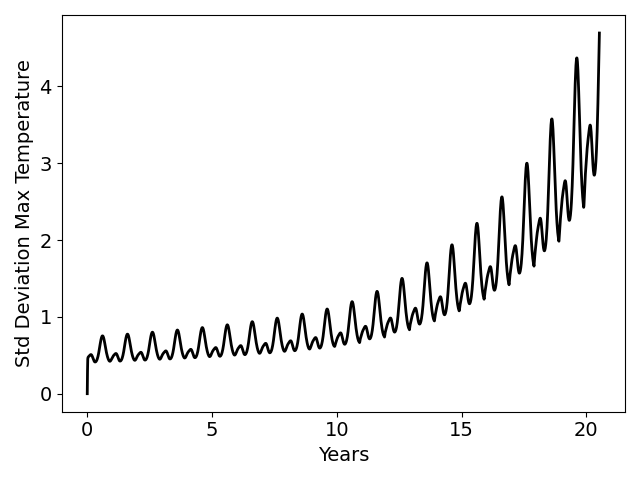}}
	\caption{Time-series of expectation and standard deviation of maximum temperature for material parameter uncertainty set $\xi_m(\omega)$ in the Texas scenario.}
	\label{fig:uq_max_1_1}
\end{figure}

\begin{figure}[H]
	\centering
	\subfloat[Expectation.]{\includegraphics[width=0.45\textwidth]{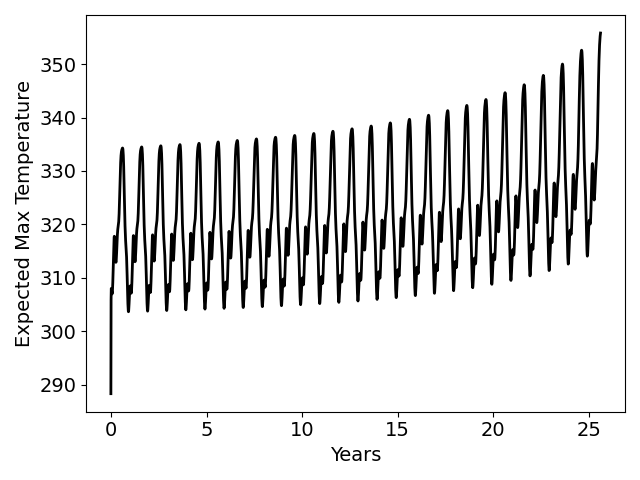}}
	\subfloat[Standard deviation.]{\includegraphics[width=0.45\textwidth]{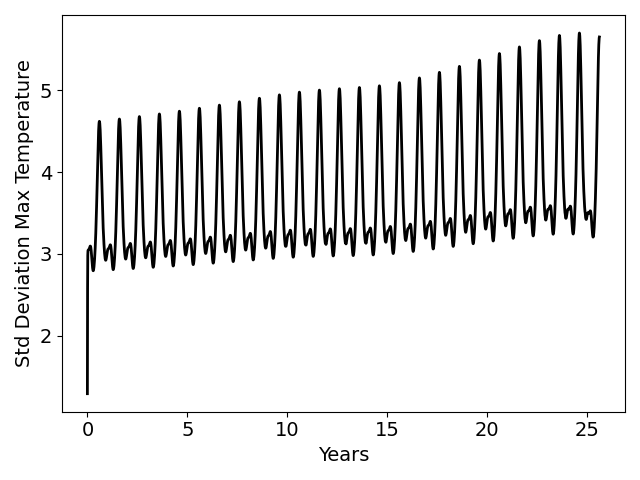}}
	\caption{Time-series of expectation and standard deviation of maximum temperature for external load uncertainty set $\xi_l(\omega)$ in the Texas scenario.}
	\label{fig:uq_max_1_2}
\end{figure}

We begin by examining the uncertainty in the material parameters from the set $\xi_m(\omega)$. We analyze the expected temperature field and its standard deviation over time, presented in 5-year increments, and show the results in Fig~\ref{fig:uq_1_1}. The results confirm that, as in the deterministic case, the maximum temperature and maximum standard deviation remain at the center of the conductor. Next, we examine the time-series evolution of maximum temperature at the center of the transmission line in Fig.\ref{fig:uq_max_1_1}. We note that the maximum temperature mean and standard deviation increase over time. We then assess the time-series evolution of maximum temperature for the set $\xi_l(\omega)$ in Fig.\ref{fig:uq_max_1_2}. We observe a steadier increase in the standard deviation with larger fluctuations.


\begin{figure}[H]
	\centering
	\subfloat[Material parameters.]{\includegraphics[width=0.45\textwidth]{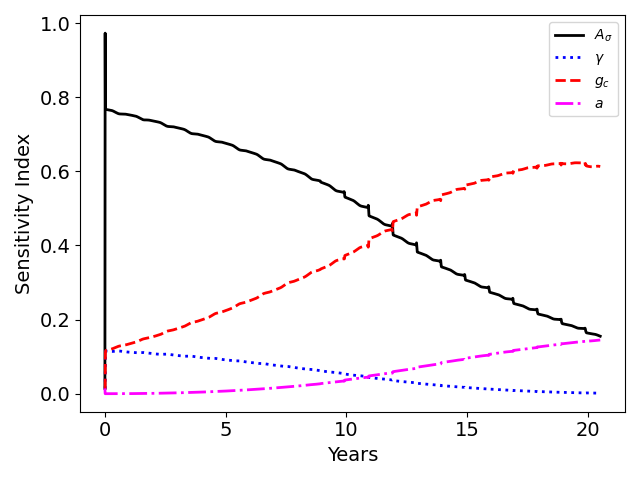}}
	\subfloat[External loading.]{\includegraphics[width=0.45\textwidth]{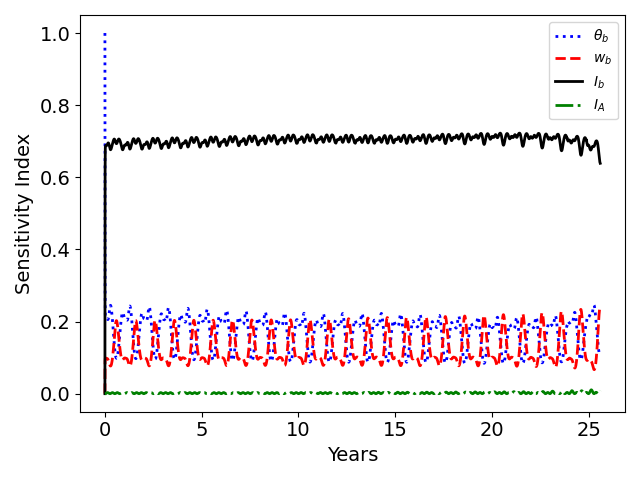}}
	\caption{First-order sensitivity index $S_i$ for the Texas scenario for parameter sets $\xi_m(\omega)$ and $\xi_l(\omega)$.}
	\label{fig:si_1_12}
\end{figure}

Using realizations from the PCM collocation points, we calculate the Sobol indices $S_i$ according to Eq.(\ref{eq:si}) and illustrate the time-series evolution of all parameters from both sets $\xi_m(\omega)$ and $\xi_l(\omega)$ in Fig\ref{fig:si_1_12}. Initially, the cross-section area parameter, which influences damage localization, is significant in initiating damage but becomes less important as the simulation progresses. Over time, $g_c$ and $a$ become more significant, as they are related to the total energy threshold for fracture and the rate of fatigue accumulation, respectively. Among the loading parameters, the current base parameter $I_b$ from the initial state emerges as the most crucial in affecting the uncertainty of $\theta_{max}$ due to its direct impact on Joule heating. Meanwhile, the wind base parameter $w_b$  and the temperature base parameter also play a significant role due to the interplay between Joule heating and convective cooling.

\begin{figure}[H]
	\subfloat[Expectation.]{\includegraphics[width=0.33\textwidth]{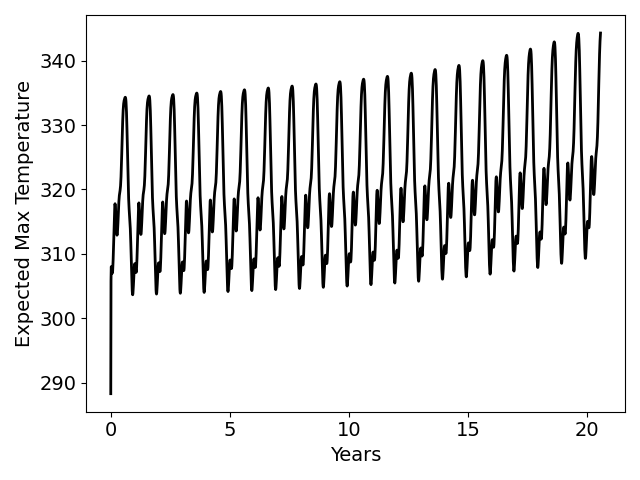}}
	\subfloat[Standard deviation.]{\includegraphics[width=0.33\textwidth]{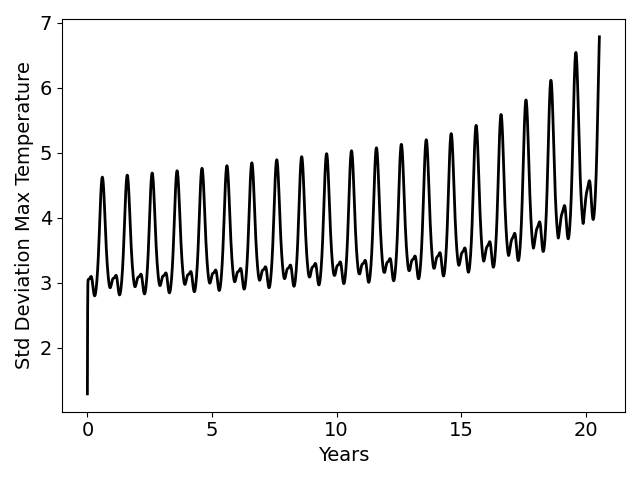}}
	\subfloat[Sensitivity index.]{\includegraphics[width=0.33\textwidth]{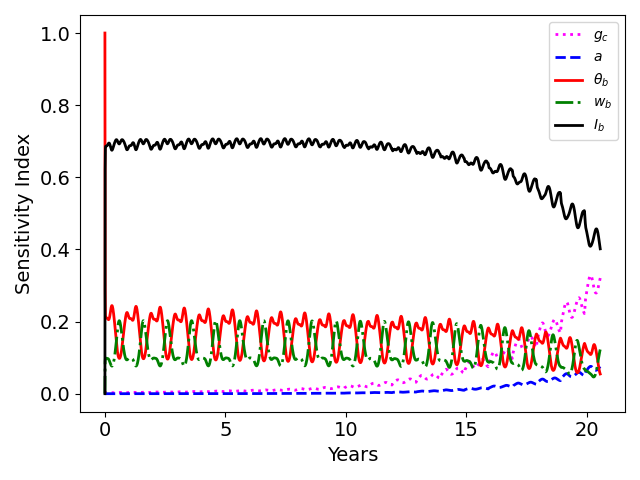}}
	\caption{Time-series of expectation, standard deviation, and sensitivity index of maximum temperature under the combined parametric space set $\xi_1(\omega)$ in the Texas scenario.}
	\label{fig:si_1_3}
\end{figure}

We combine the two most influential material parameters with the three most influential loading condition parameters to create a new set, $\xi_1 (\omega) = {g_c(\omega), a(\omega),\theta_b(\omega),w_b(\omega),I_b(\omega)}$, and conduct a final round of sensitivity analysis for Texas. The results in Fig.~\ref{fig:si_1_3} show that initially, the loading conditions have a more substantial impact than material parameters, but over time, as aging effects become more significant, the relative importance of material parameters increases.

\subsubsection{California, Michigan, and Florida}
In California, Michigan, and Florida scenarios, we combine the set of random parameters $\xi (\omega) = \{g_c(\omega), a(\omega),\theta_b(\omega),w_b(\omega),I_b(\omega)\}$, and perform a final round of SA.

\begin{figure}[H]
	\subfloat[Expectation.]{\includegraphics[width=0.33\textwidth]{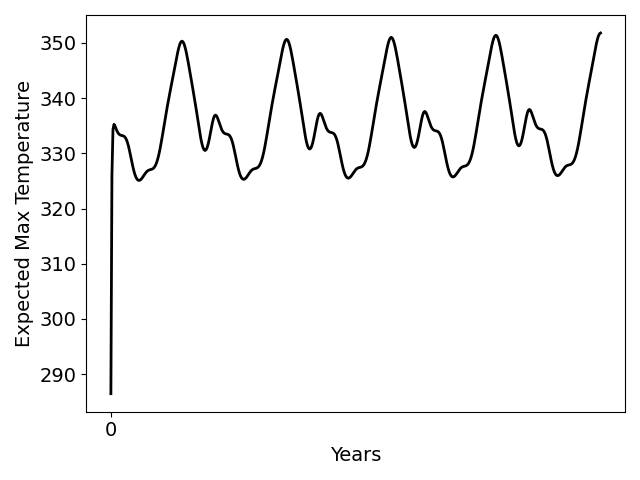}}
	\subfloat[Standard deviation.]{\includegraphics[width=0.33\textwidth]{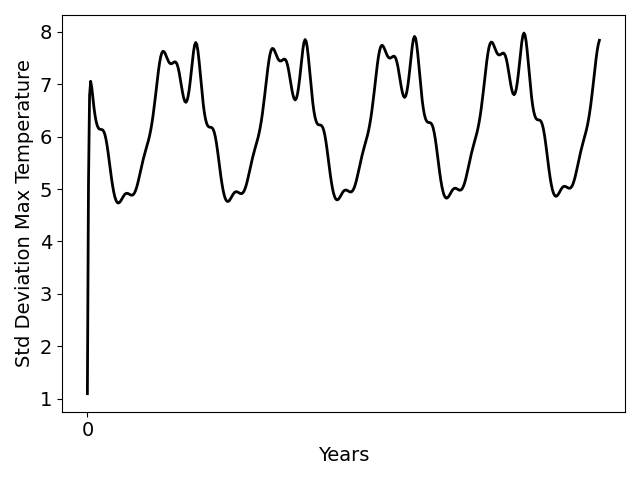}}
	\subfloat[Sensitivity index.]{\includegraphics[width=0.33\textwidth]{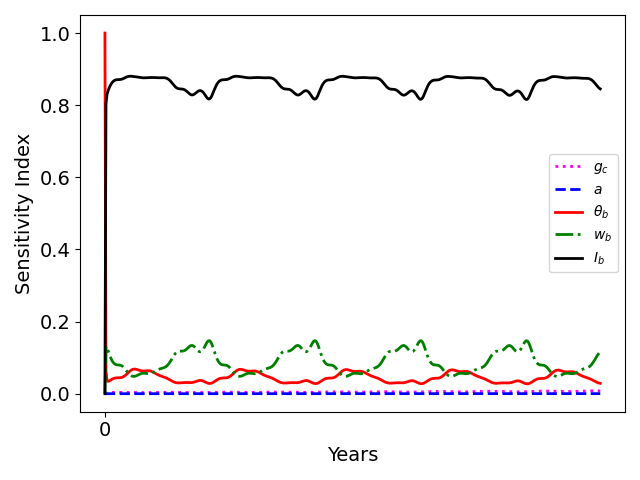}}
	\caption{Time-series of expectation, standard deviation, and sensitivity index of maximum temperature under the combined parametric space in the California scenario.}
	\label{fig:si_2}
\end{figure}

In the California scenario, when the time series for maximum temperature is truncated at the earliest occurrence of failure, the analysis reveals that failures happen within a span shorter than five years. However, the pattern in Fig.~\ref{fig:si_2} aligns closely with observations from the Texas Scenario. The consistency in the pattern highlights the dominant influence of the current base parameter, $I_b$, which is intensified by low wind speeds that reduce the convective cooling. As a result, the conductor temperature rapidly increases to critical thresholds, driven by the Joules heating in the absence of effective cooling. 

\begin{figure}[H]
	\subfloat[Expectation.]{\includegraphics[width=0.33\textwidth]{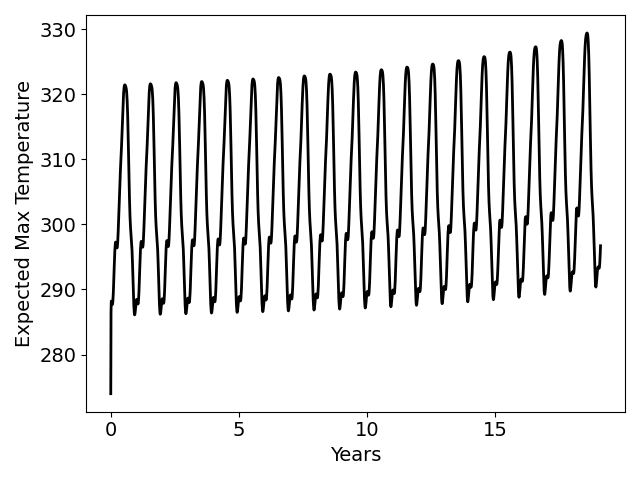}}
	\subfloat[Standard deviation.]{\includegraphics[width=0.33\textwidth]{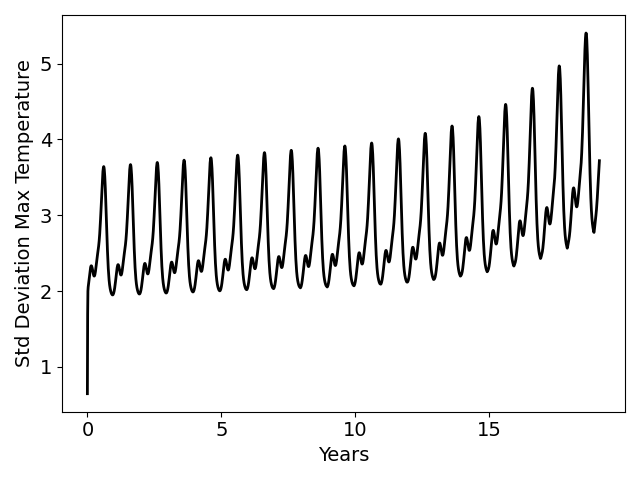}}
	\subfloat[Sensitivity index.]{\includegraphics[width=0.33\textwidth]{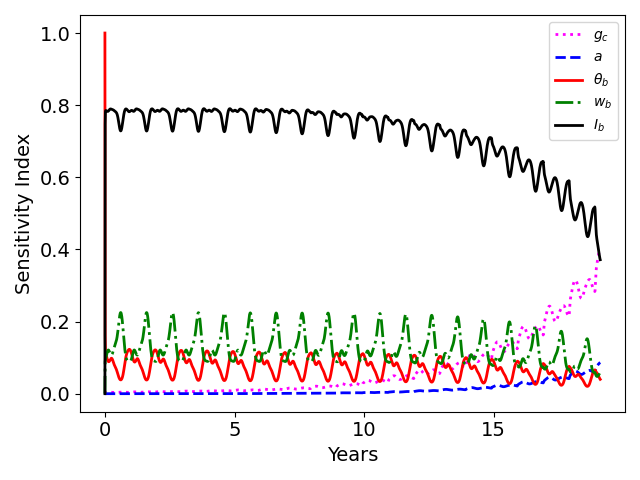}}
	\caption{Time-series of expectation, standard deviation, and sensitivity index of maximum temperature under the combined parametric space in the Michigan scenario.}
	\label{fig:si_3}
\end{figure}

In the Michigan scenario, shown in Fig.~\ref{fig:si_3}, the current base parameter, $I_b$, initially has the most significant impact on the failure of the transmission line. Over time, however, $g_c$ becomes increasingly influential, eventually surpassing the impact of $I_b$. Historical data from Michigan shows that the state experiences higher wind speeds compared to the other states analyzed. Although the wind has a convective effect, the impact of increased mechanical load on the cable is more significant accelerating the aging of the transmission line.

\begin{figure}[H]
	\subfloat[Expectation.]{\includegraphics[width=0.33\textwidth]{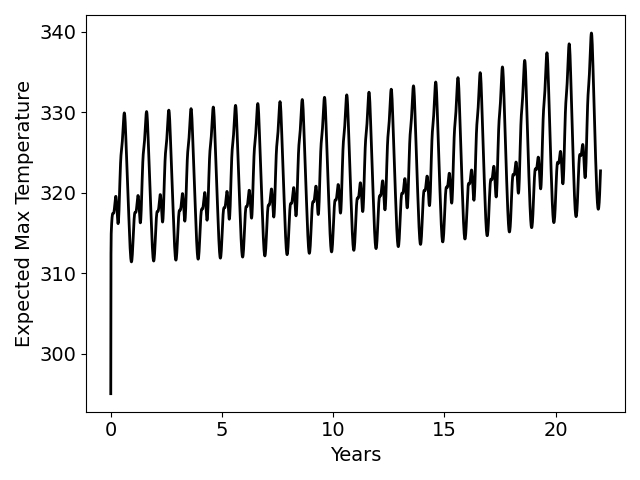}}
	\subfloat[Standard deviation.]{\includegraphics[width=0.33\textwidth]{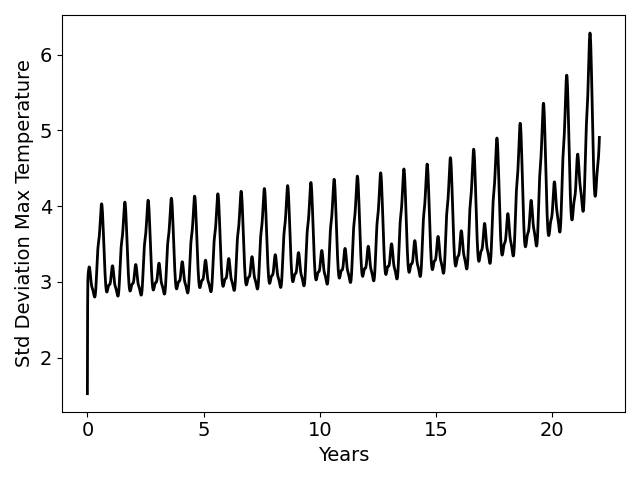}}
	\subfloat[Sensitivity index.]{\includegraphics[width=0.33\textwidth]{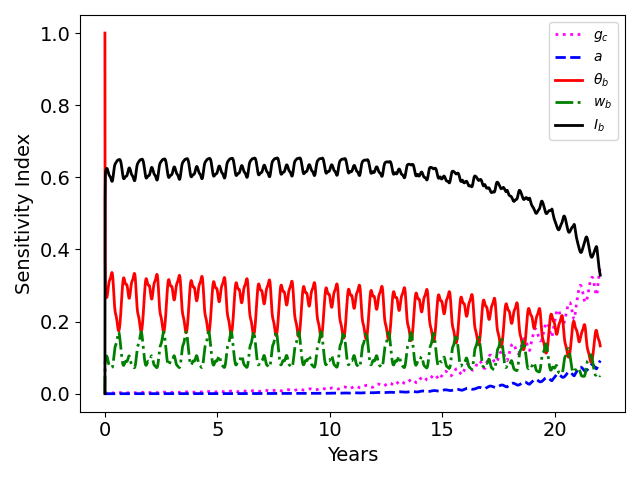}}
	\caption{Time-series of expectation, standard deviation, and sensitivity index of maximum temperature under the combined parametric space in the Florida scenario.}
	\label{fig:si_4}
\end{figure}

In the Florida Scenario, as illustrated in Fig.~\ref{fig:si_4}, the current base parameter, $I_b$, initially shows the most significant influence on the transmission line's failure. However, over time, the fracture energy parameter, $g_c$, becomes more dominant. Historical data considered in this study shows that Florida experiences high temperatures throughout the year. Despite higher wind speeds, the effectiveness of convective cooling is diminished due to the high ambient temperatures.

\subsubsection{Probability of Failure}
\begin{figure}[H]
	\centering
	\subfloat[Probability of Failure]{\includegraphics[width=0.45\textwidth]{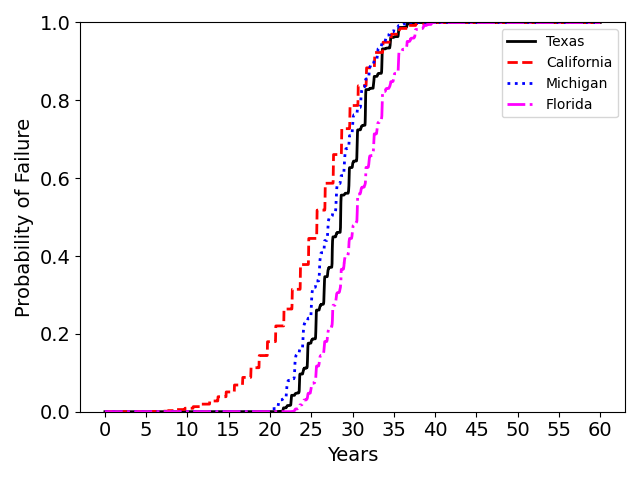}}
	\caption{Probability of failure for Texas, California, Michigan, and Florida.}
	\label{fig:pf_all}
\end{figure}

In this section, we calculate the expected value of the Bernoulli variable $h$, as specified in Eq.(\ref{eq:h}), using the PCM framework with $n=5$ points based on the multi-dimensional expectation formula in Eq.(\ref{eq:pcm_multi}). The results are displayed in Fig~\ref{fig:pf_all} for the reference mean parameter values from the parameter sets $\xi_1(\omega)$, $\xi_2(\omega)$, $\xi_3(\omega)$, and $\xi_4(\omega)$. Under the consideration of four specific states, the probability of failure curve for California initially shows the early higher possibility of occurrence of failure. In contrast, the Florida scenario shows the failure curve on the right shift indicating a lower possibility of failure among the four states.
Finally, we analyze the impact of varying mean values of five influential parameters $\xi_1 (\omega) = {g_c(\omega), a(\omega),\theta_b(\omega),w_b(\omega),I_b(\omega)}$ under three different levels of $A_\sigma$ in the time-series probability of failure ($p_f$) for each scenario. The comparison includes a baseline scenario representing a transmission line with minimal initial damage. Unlike previous sensitivity analyses of parameters, which provided insights into levels of importance relative to their uncertainty, this section focuses on the effects of changing the baseline values of influential parameters on the $p_f$ curves under minimal, moderate, and severe initial damage.

\begin{figure}[H]
    \centering
    \subfloat[Texas]
    {\includegraphics[width=0.45\textwidth]{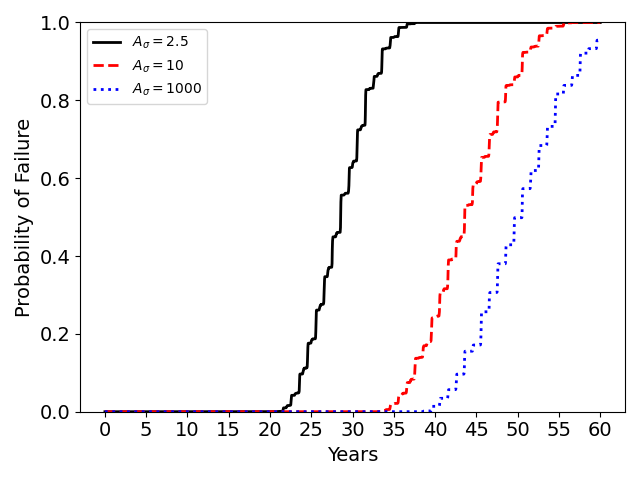}}
    \subfloat[California]{\includegraphics[width=0.45\textwidth]{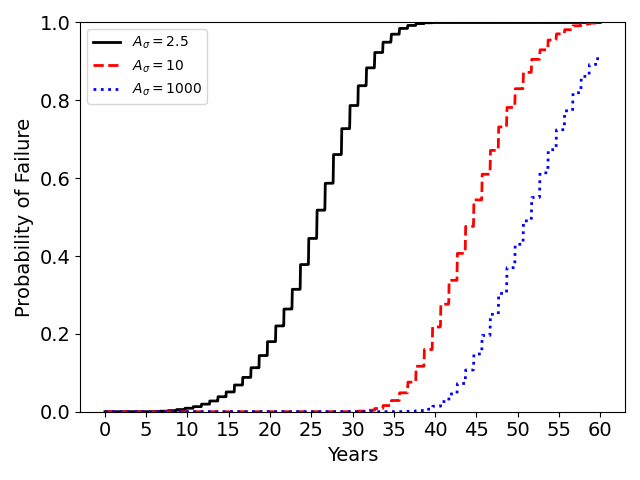}}
    \hfill 
    \subfloat[Michigan]
    {\includegraphics[width=0.45\textwidth]{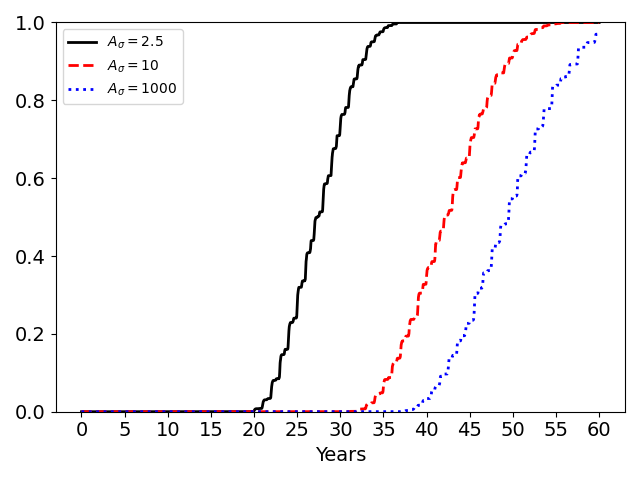}}
    \subfloat[Florida]
    {\includegraphics[width=0.45\textwidth]{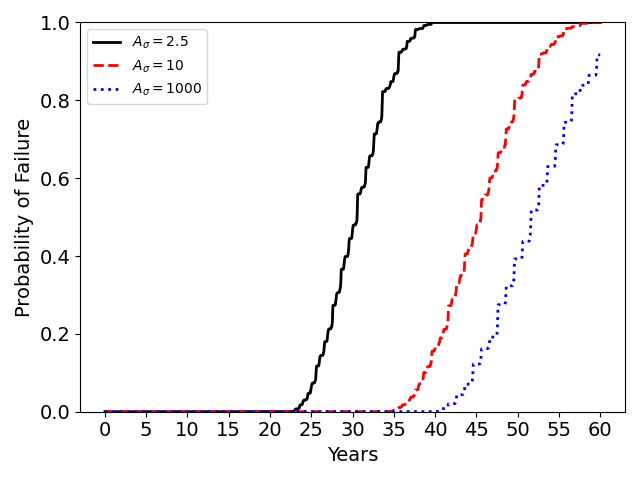}}
    \caption{Probability of Failure of transmission line for different values of initial damage (shown in legends).}
    \label{fig:individual Failure}
\end{figure}

In each scenario, we observe three distinct curve shapes representing different levels of damage, as depicted in Fig.\ref{fig:individual Failure}. The analysis shows that in the presence of moderate initial damage, the chance of failure increases from $20 \%$ to $60 \%$ after 45 years in the Texas, California, and Michigan scenarios. In Florida, the probability of failure increases from $20 \%$ to $60 \%$ after 45 years. In the presence of severe initial damage, the probability of failure reaches $100\%$ well before the failure initiates in a minimal damage case. All the states show the life span of the transmission lines to be around 35 to 40 years under severe damage. This indicates how detrimental the effect of initial damage could be on the life span of the transmission line. Among all the four scenarios, the chance of failure initiated early in the case of California as it has low wind and high ambient temperature reducing the impact of convective cooling. Also, the plots reveal that the probability of failure is around $90\%$ in 60 years indicating there is a $10\%$ chance the transmission line will last more than 60 years. However, in the case of Michigan, as it has higher fluctuations in wind and temperature, the chance of failure is maximum at 60 years in comparison to the three states.

\subsection{Convergence Analysis}
\begin{figure}[H]
	\centering
	\subfloat[PCM.]{\includegraphics[width=0.45\textwidth]{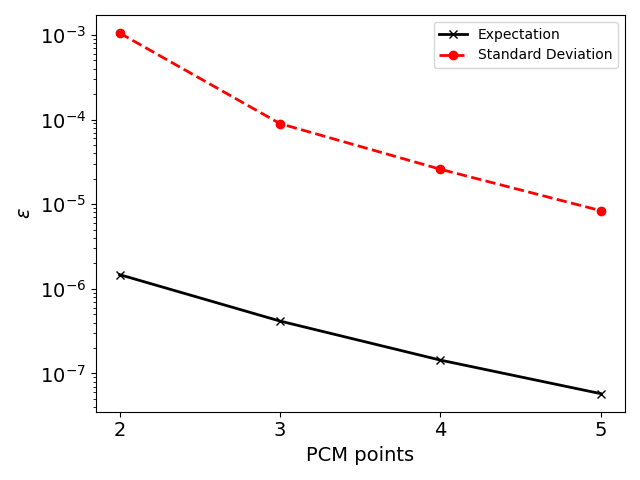}}
      \subfloat[MC.]{\includegraphics[width=0.45\textwidth]{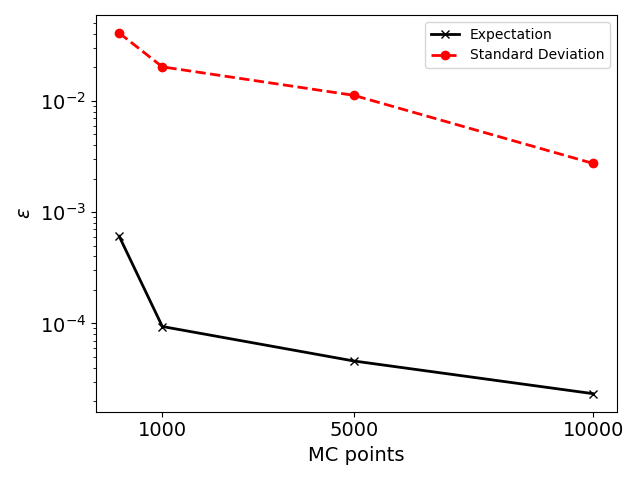}}
	\caption{Convergence Plots.}
	\label{fig:convergence}
\end{figure}

The results provide a compelling perspective for analyzing the reliability of transmission lines. However, it is important to verify the consistency of our PCM results through convergence analysis. As the model involves the coupling of several governing equations, obtaining an analytical solution is not feasible. Therefore, we depend on a refined PCM solution to serve as our reference.  
\begin{equation}
\epsilon = \frac{\Vert \theta -\theta_{ref}\Vert_2 }{\Vert \theta_{ref} \Vert_2}.
\end{equation}
In our convergence study, we focus on the effects of the most influential parameter $I_b$ simplifying the Probabilistic Collocation Method (PCM) to a 1-D problem, enabling us to achieve accurate integration using 100 collocation points as the reference solution. We then compare the performance of lower-order PCM solutions to Monte Carlo simulations, analyzing the relative errors in the norm of the temperature field solution at a specific time. The comparison is graphically represented in  Fig~\ref{fig:convergence} to highlight the accuracy of PCM relative to Monte Carlo methods.
We calculated the relative errors until the earliest failure considering the Texas Scenario as a reference. The Probabilistic Collocation Method achieved significantly lower errors, approximately two orders of magnitude less, compared to 10,000 Monte Carlo simulations using just 5 collocation points. Using this strategy is particularly effective as we move to higher-dimensional spaces, where MC methods would require more realizations to represent the space accurately.

\section{Conclusion}
We developed a thermo-electro-mechanical model to predict the failure of transmission lines across four specific states in the US: Texas, California, Michigan, and Florida. We considered the historical wind and temperature data of each state. The mechanical model evaluates material damage and aging from long-term fatigue. It also integrates a thermal model that includes a heat transfer equation, accounting for Joule heating due to electric currents and convective cooling due to the wind. Additionally, the electrical model addresses how accumulated damage and temperature-induced resistance cause voltage drops along the transmission line. Overall, the model acts as a positive feedback loop, where initially damaged transmission lines deteriorate to the point where the material reaches its annealing temperature, ultimately leading to failure. 

We used the finite element method to solve a set of quasi-static equations. We studied four different states of the US: Texas, California, Michigan, and Florida to understand the long-term behavior of the transmission lines in the presence of minimal, moderate, and severe damage. Each scenario was analyzed deterministically to track the evolution of the maximum conductor temperature under varying wind, temperature, and current loading. We used discrete Fourier analysis to obtain the cyclic loading condition from the discrete historical data of wind and temperature. To reduce the complexity, we parameterized the current loading. Subsequently, we utilized the Probabilistic Collocation Method (PCM) for uncertainty quantification (UQ), sensitivity analysis (SA), and probability of failure assessments.

The deterministic solution revealed how the temperature in the conductor under cyclic wind, air temperature, and current loading in the presence of initial damage affects the overall longevity of the transmission lines. Even the moderate initial damage significantly reduced the material's life expectancy. 
In the global sensitivity analysis using PCM, we identified the electric current, $I_b$ initially has a dominating effect which was later surpassed by fracture energy, $g_c$ over time.
Additionally, with the application of PCM, by defining the limit state function as a Bernoulli random variable, we studied the probability of failure for each specific scenario, varying the most influential parameters for three damage cases. The analysis revealed, that the presence of initial damage significantly increased the chance of failure.

We acknowledge that the current model incorporates simplifying assumptions that may be further explored in the future. For instance, replacing the quasi-static equations with a transient analysis could provide more dynamic insights into the system's response. Furthermore, while the cable mechanics are currently simplified to one dimension, an extension to model the catenary shape could enhance the realism of our simulations. Even the integration of real current data for specific states along with the effect of solar radiation and solar heat gain could enhance the predictive accuracy of failure probabilities. All these modifications shift the dependency from abstract stochastic process modeling to a more robust physics-based approach, significantly improving the reliability assessments of the transmission line.

\bibliographystyle{ieeetr}
\bibliography{references}
\end{document}